\documentclass[12pt,a4paper]{article}

\usepackage[utf8]{inputenc}
\usepackage[english]{babel}
\usepackage{amsmath,amssymb,amsfonts,amsthm,graphicx}
\usepackage{hyperref}
\usepackage{cases}

\begin{document}

\newcommand{\C}{\mathbb {C}}
\newcommand{\N}{\mathbb{N}}
\newcommand{\R}{\mathbb{R}}
\newcommand{\tf}{\mathcal{F}}

\renewcommand{\theenumi}{\roman{enumi}}
\renewcommand{\labelenumi}{\theenumi)}

\swapnumbers
\newtheorem{def1}{Definition}[section]
\newtheorem{thm1}[def1]{Theorem}
\newtheorem{thm2}[def1]{Theorem}
\newtheorem{thm3}{Theorem}[section]
\newtheorem{rem2}[thm3]{Remarks}
\newtheorem{thm4}[thm3]{Theorem}
\newtheorem{rem3}[thm3]{Remark}
\newtheorem{prop1}{Proposition}[section]
\newtheorem{thm5}[prop1]{Theorem}
\newtheorem{rem4}[prop1]{Remark}
\newtheorem{rem1}[prop1]{Remark}
\newtheorem{cor2}[prop1]{Theorem}
\newtheorem{prop2}{Proposition}[section]
\newtheorem{prop3}[prop2]{Proposition}
\newtheorem{prop4}[prop2]{Proposition}
\newtheorem{cor1}[prop2]{Proposition}

\title{Explicit error estimates for the stationary phase method II:\\ Interaction of amplitude singularities with stationary points}

\author{Felix Ali Mehmeti \footnote{Université de Valenciennes et du Hainaut-Cambrésis, LAMAV, FR CNRS 2956, Le Mont Houy, 59313 Valenciennes Cedex 9, France. Email: felix.ali-mehmeti@univ-valenciennes.fr} , Florent Dewez \footnote{Université Lille 1, Laboratoire Paul Painlevé, CNRS U.M.R 8524, 59655 Villeneuve d'Ascq Cedex, France. Email: florent.dewez@math.univ-lille1.fr}}

\date{}


\maketitle

\begin{abstract}
	In this paper, we improve slightly Erdélyi's version of the stationary phase method by replacing the employed smooth cut-off function by a characteristic function, leading to more precise remainder estimates. We exploit this refinement to study the time-asymptotic behaviour of the solution of the free Schrödinger equation on the line, where the Fourier transform of the initial data is compactly supported and has a singularity. We obtain uniform estimates of the solution in space-time regions which are asymptotically larger than any space-time cones. Moreover we expand the solution with respect to time on the boundaries of the above regions, showing the optimality of the decay rates of the estimates.
\end{abstract}

\vspace{0.3cm}

\noindent \textbf{Mathematics Subject Classification (2010).} Primary 41A80; Secondary 41A60, 35B40, 35B30, 35Q41.

\noindent \textbf{Keywords.} Asymptotic expansion, stationary phase method, error estimate, Schrö\-din\-ger equation, $L^{\infty}$-time decay, singular frequency, space-time cone.

\setcounter{section}{-1}
\section{Introduction}

\hspace{2.5ex} The asymptotic behaviour of oscillatory integrals with respect to a large parameter, sometimes used to study long-time asymptotics for solutions of dispersive equations, can often be described using the stationary phase method. A theorem of A. Erdélyi \cite[section 2.8]{erdelyi} permits to treat oscillatory integrals with singular amplitudes and furnishes asymptotic expansions with explicit remainder estimates. In part I of this article \cite{article1}, we provided a complete proof of an improved version and we applied this result to the solution of the free Schrödinger equation on the line for initial conditions with singular Fourier transform with support in a compact interval. We obtained asymptotic expansions to one term with respect to time in certain space-time cones. We observed a blow-up of the remainder estimates when the boundaries of the cones tend to certain critical directions related to the compact interval. This phenomenon limits the regions where we can estimate uniformly the solution. In this paper, we slightly improve the above stationary phase method in order to enlarge the regions of uniform estimates. For this purpose, we replaced the smooth cut-off function used in \cite{article1} and \cite{erdelyi}, which contributes partially to the blow-up, by a characteristic function. This permits to make explicit the blow-up. We exploit the obtained refinement to provide estimates of the solution of the free Schrödinger equation in space-time regions which are asymptotically larger than the cones considered in \cite{article1}, and to expand the solution with respect to time on curves in space-time, where decay and blow-up balance out.

Let us introduce the free Schrödinger equation
\begin{equation*}
	(S) \qquad \left\{ \begin{array}{rl}
			& \big[ i \partial_t + \partial_x^2 \big] u(t,x) = 0 \\
			& \vspace{-0.3cm} \\
			& u(0,x) = u_0(x)
	\end{array} \right. \; ,
\end{equation*}
for $t > 0$ and $x \in \R$. If we suppose $u_0 \in L^{1}(\R)$ then
\begin{equation*}
	\left\| u(t,.) \right\|_{L^{\infty}(\R)} \leqslant \frac{\| u_0 \|_{L^1(\R)}}{2\sqrt{\pi}} \, t^{-\frac{1}{2}} \; ,
\end{equation*}
see for example \cite[p.60]{reed-simon}. If it is assumed that $u_0 \in L^2(\R)$ then we have by Strichartz' estimate (\cite{strichartz}, see also \cite{banica})
\begin{equation*}
	\left\| u(t,.) \right\|_{L^{\infty}(\R)} \leqslant C \, \| u_0 \|_{L^2(\R)} \, t^{-\frac{1}{4}} \; ,
\end{equation*}
for some suitable constant $C>0$. Now, consider for example initial conditions $u_0$ satisfying
\begin{equation} \label{intro}
	\forall \, p \in \R \qquad	\tf u_0 (p) = p^{\mu-1} (1-p) \, \chi_{[0,1]}(p) \; ,
\end{equation}
with $\mu \in (0,1)$; here $\tf u_0$ refers to the Fourier transform of $u_0$ and $\chi_{[0,1]}$ is the characteristic function of $[0,1]$. Under this assumption, $u_0$ is a smooth function which never belongs to $L^1(\R)$ and belongs to $L^2(\R)$ if and only if $\mu \in \big(\frac{1}{2}, 1\big)$. In \cite{article1}, we showed that there exists a constant $c(u_0,\varepsilon_1,\varepsilon_2) \geqslant 0$ such that for all $(t,x)$ lying in the space-time cone  $\mathfrak{C}_{\varepsilon_1,\varepsilon_2}(0,1) := \big\{ (t,x) \, \big| \, t>0 \, , \, \varepsilon_1 \leqslant \frac{x}{2t} \leqslant 1-\varepsilon_2 \big\}$, with fixed $\varepsilon_1, \varepsilon_2 > 0$ satisfying $\varepsilon_1 < 1 - \varepsilon_2$, there exist $H(t,x,u_0), K_{\mu}(t,x,u_0) \in \C$ such that
\begin{equation} \label{exple1}
	\left| u(t,x) - H(t,x,u_0) \, t^{-\frac{1}{2}} - K_{\mu}(t,x,u_0) \, t^{-\mu} \right| \leqslant c(u_0,\varepsilon_1,\varepsilon_2) \, \big(t^{-1} + t^{-\delta} \big) \; ,
\end{equation}
where $\delta \in \big( \max\{\mu,\frac{1}{2}\},1\big)$. This result shows that the preponderant time-asymptotic behaviour of the solution depends on the value of $\mu$ and we derive the optimal decay rate in the cone, which is either $t^{-\frac{1}{2}}$ or $t^{-\mu}$. Moreover on the space-time direction given by the singular frequency, as for example $\frac{x}{2t} = 0$ in \eqref{intro}, we proved that the decay rate of the solution is $t^{-\frac{\mu}{2}}$. This slow decay is the consequence of the interaction between the stationary point and the singularity. Nevertheless the constant $c(u_0,\varepsilon_1,\varepsilon_2)$ which appears in \eqref{exple1} blows up when $\varepsilon_1$ or $\varepsilon_2$ tends to $0$, and so it seems to be impossible to deduce from these results a uniform $L^{\infty}$-time decay rate. Our motivation in this paper is to enlarge the controllable space-time regions and to understand better the transition between the behaviour in the cone and the behaviour on the critical direction.

In section 1, we provide a stationary phase method based on \cite{article1} for oscillatory integrals of the type
\begin{equation*}
	\int_{p_1}^{p_2} U(p) \, e^{i \omega \psi(p) } \, dp \; ,
\end{equation*}
where the endpoints of the interval are integrable singularities of the amplitude $U$ and stationary points of real order of the phase $\psi$. Here we replace the smooth cut-off function used in \cite{article1} and \cite{erdelyi} to separate the different behaviours by a characteristic function. For this purpose, we consider a fixed cutting-point $q \in (p_1,p_2)$. After the splitting of the integral, we carry out substitutions to simplify the phases and we integrate by parts to create asymptotic expansions to one term at $p_1$ and $p_2$, following the lines of \cite[section 1]{article1}. Here the characteristic function produces terms related to the cutting-point. Due to the different substitutions employed in the two integrals, these terms are not the same. However by expanding them with respect to the parameter $\omega$, we observe  that the first terms cancel out but not the remainders. Finally, we estimate these remainders related to the cutting point as well as the remainders coming from the expansions at $p_1$ and $p_2$.\\
Assuming that the endpoint $p_j$, for a certain $j \in \{1,2\}$, is a regular point of the amplitude, the previous result furnishes only the same decay rate for the first term and the remainder related to $p_j$. Consequently, the refined estimate of the remainder given in \cite[Theorem 1.7]{article1}, which provides a better decay rate, is employed. Moreover, we note that the rest related to the point $q$ is always negligible as compared with the first term.

In section 2, we apply the above method to an oscillatory integral where the phase is a polynomial function of degree 2 and where the amplitude has a unique singularity at the left endpoint of the interval, in view of an application to the free Schrödinger equation. Thanks to the explicitness of the phase function and the preciseness of the previous results, we furnish an asymptotic expansion together with remainder estimates which depend explicitly on the distance between the two endpoints of the interval. For fixed $\omega$, a blow-up of the first terms and the remainder estimates occurs when the stationary point approaches the singularity. On the other hand, if we suppose that the distance between the stationary point and the singularity is lower bounded by a positive constant, then the first terms and the remainder both tend to $0$ at a different rate when the large parameter $\omega$ tends to infinity. This can be exploited by assuming that the stationary point tends to the singularity as the large parameter tends to infinity. Under this hypothesis, we provide asymptotic expansions of the oscillatory integral on curves in the space of the parameters. These curves are indexed by some positive $\varepsilon>0$ describing the convergence speed of the stationary point towards the singularity, which is chosen in such a way that the decay rate of the remainder is always larger than the decay rate of the first terms.

In section 3, we consider the Fourier solution formula of the free Schrödinger equation on the line with initial conditions in a compact frequency band with a singular frequency at one of the endpoints. We employ the precise remainder estimates found in section 2 to derive an estimate of the solution in regions where the distance between the stationary point and the singularity is bounded from below by a small negative power of the time. We point out that these regions are asymptotically larger than any space-time cone and so the space-time regions with uniform estimates are enlarged, as compared with \cite{article1}. To show the optimality of these estimates, we give asymptotic expansions of the solution formula with respect to time on space-time curves, which are the boundaries of the above regions.

In section 4, we remind some technical results established in \cite{article1} which have been employed in section 1.

Applied to example \eqref{intro}, our results yield either
\begin{equation*} \label{ASYMPT1}
		\Big| u(t,x) - H(t,u_0) \, t^{-\frac{1}{2} + \varepsilon(1-\mu)} \, \Big| \leqslant R_{\mu > 1/2}(u_0,\varepsilon) \, t^{-\mu + \varepsilon \mu} \; ,
\end{equation*}
if $\mu > \frac{1}{2}$, or
\begin{equation*} \label{ASYMPT2}
		\Big| u(t,x) - K_{\mu}(t,u_0) \, t^{-\mu + \varepsilon \mu} \, \Big| \leqslant R_{\mu < 1/2}(u_0,\varepsilon) \, t^{-\frac{1}{2} + \varepsilon(1-\mu)} \; ,
\end{equation*}
if $\mu < \frac{1}{2}$, for all $(t,x)$ belonging to the curve $\mathfrak{G}_{\varepsilon} := \big\{ (t,x) \, \big| \, t > 1 \, , \, \frac{x}{2t} = t^{-\varepsilon} \big\}$, with fixed $\varepsilon \in \big(0, \frac{1}{2} \big)$ and
\begin{align*}
	& \bullet \quad H(t,u_0) := \frac{1}{2 \sqrt{\pi}} \, e^{-i \frac{\pi}{4}} \,  e^{i t^{1-\varepsilon}}  \left(1-t^{-\varepsilon} \right) \; , \\
	& \bullet \quad K_{\mu}(t,u_0) := \frac{\Gamma(\mu)}{2^{\mu+1} \pi} \, e^{i \frac{\pi \mu}{2}} \; .
\end{align*}
The case $\mu = \frac{1}{2}$ is also studied in this paper. We observe that the parameters $\mu$ and $\varepsilon$ affect the decay rate of the solution. In \cite{article1}, we have already noticed that the strength of the singularity has an effect on the decay. Here we remark that the diminution of the positive quantity $\frac{1}{2} - \varepsilon$ increases the convergence speed of the stationary point towards the singularity and thus the influence of the blow-up. We emphasize this phenomenon by noting that the constants $R_{\mu > 1/2}(u_0,\varepsilon), R_{\mu < 1/2}(u_0,\varepsilon)$ blow up when the parameter $\varepsilon$ reaches the critical value $\frac{1}{2}$. This illustrates the transition between blow-up and decay.\\
We also obtain estimates of the solution in the space-time region $\mathfrak{R}_{\varepsilon} := \big\{ (t,x) \, \big| \, t^{-\varepsilon} \leqslant \frac{x}{2t} < 1 \, , \, t > 1 \big\}$, namely either
\begin{equation} \label{ESTGLO1}
	\big| u(t,x) \big| \leqslant C_{\mu > 1/2}(u_0,\varepsilon) \, t^{-\frac{1}{2} + \varepsilon(1-\mu)} \; ,
\end{equation}
if $\mu > \frac{1}{2}$, or
\begin{equation} \label{ESTGLO2}
	\big| u(t,x) \big| \leqslant C_{\mu < 1/2}(u_0,\varepsilon) \, t^{-\frac{1}{2} + \varepsilon(1-\mu)} \; ,
\end{equation}
if $\mu < \frac{1}{2}$, for a fixed $\varepsilon \in \big(0,\frac{1}{2} \big)$. We observe that the decay rates are attained on the curve $\mathfrak{G}_{\varepsilon}$, the left boundary of $\mathfrak{R}_{\varepsilon}$, which shows the optimality of the rates on $\mathfrak{R}_{\varepsilon}$. Moreover the decay rates in \eqref{ESTGLO1} and \eqref{ESTGLO2} tend to $t^{-\frac{\mu}{2}}$ when $\varepsilon$ tends to $\frac{1}{2}$, which is the decay rate of the solution on the direction given by the singular frequency according to \cite{article1}. But the constants $C_{\mu > 1/2}(u_0,\varepsilon), C_{\mu < 1/2}(u_0,\varepsilon)$ blow up as $\varepsilon$ tends to $\frac{1}{2}$, for the same reason as above.

These results lead us to the conjecture that a global optimal estimate
\begin{equation*}
	\left\| u(t,.) \right\|_{L^{\infty}(\R)} \leqslant C(u_0) \, t^{-\frac{\mu}{2}}
\end{equation*}
holds. This decay rate would be optimal according to the previous results. Further a generalization of the concentration phenomenon around critical directions caused by singular frequencies seems to be possible, as for example in the context of other propagation phenomena or in the context of other geometries.

This paper, together with its part I, see \cite{article1}, have been inspired by \cite{fam}, where the authors calculated asymptotic expansions to one term with respect to time of solutions of the Klein-Gordon equation on a star shaped network. The coefficients are constant but different on the semi-infinite branches. The authors determined the influence of these coefficients on the expansion.

\cite{fam} shows the way to obtain this type of results for the Schrödinger equation on domains with canonical geometry and canonical potential permitting sufficiently explicit spectral theoretic solution formulas.

Our refined version of Erdélyi's expansion theorem could help to improve the comprehension of a difficulty described in \cite{fam}: when the frequency band approaches certain critical values coming from potential steps, the asymptotic expansion to one term degenerates. These critical values will play a similar role as the singularities of the Fourier transform of these initial conditions in the present paper.

In \cite{fam94}, a global (probably not optimal) $L^{\infty}$-time decay estimate for the Klein-Gordon equation on $\R$ with potential steps has been proved using the van der Corput inequality, spectral theory and the methods of \cite{msw}.

As further perspectives for future research, let us mention the global (and optimal) time decay estimate in the setting of this paper and the same issue for the Klein-Gordon equation on the star shaped network \cite{fam}.

Further one could review existing results on the Schrödinger equation with localized potential, e.g. \cite{weder}, \cite{goldbergschlag} and \cite{fam0}, and non linear equations, e.g. \cite{strunk}, to check the possible use of refined estimates of the free equation. An interesting issue could be to find optimal decay conditions of the potential.

The time-decay rate of the free Schrödinger equation is considered in \cite{cazenave1998} and \cite{cazenave2010}. In \cite{cazenave1998}, singular initial conditions are constructed to derive the exact $L^p$-time decay rates of the solution, which are slower than the classical results for regular initial conditions. In \cite{cazenave2010}, the authors construct initial conditions in Sobolev spaces (based on the Gaussian function), and they show that the related solutions have no definite $L^p$-time decay rates, nor coefficients, even though upper estimates for the decay rates are established.\\
The papers \cite{cazenave1998} and \cite{cazenave2010} use special formulas for functions and their Fourier transforms, which are themselves based on complex analysis. In our results, we furnish slower decay rates by considering initial conditions with singular Fourier transforms. Here, complex analysis is directly applied to the solution formula of the equation, which permits to obtain results for a whole class of functions. The method seems to be more flexible.

One can mention the stationary phase method of \cite[chapter 7]{hormander}. It is assumed that the amplitude $U$ and the phase have a certain regularity on $\R^d$ (for $d \geqslant 1$) and that $U$ has a compact support. Asymptotic expansions of the oscillatory integral are given by using Taylor's formula of the phase, where the stationary point is supposed non-degenerate. However, stronger hypotheses concerning the phase are required in order to bound uniformly the remainder by a constant, which is not explicit.

In \cite{liess}, the authors are interested in the decay of Fourier transforms on singular surfaces. They cite Erdélyi's result but they employ only the possibility to have stationary point of general integer order.

One-dimensional Schrödinger equations with singular coefficients are studied in \cite{banica}. Dispersion inequalities and Strichartz-type estimates are furnished and we observe that the singularities do not influence the results as compared with the regular case.\\

\noindent \textbf{Acknowledgements:}\\
The authors thank R. Haller-Dintelmann and V. Régnier for discussing the possibility to use characteristic functions in the context of the stationary phase method.


\section{Lossless error estimates}

\hspace{2.5ex} We introduce the two assumptions related to the amplitude and to the phase.

Let $p_1,p_2$ be two real numbers such that $-\infty < p_1 < p_2 < +\infty$.\\ \\
\textbf{Assumption (P$_{\rho_1,\rho_2,N}$).} For $\rho_1, \rho_2 \geqslant 1$, let $\psi \in \mathcal{C}^1\big([p_1,p_2], \R\big)$ be a function satisfying
	\begin{equation*}
		\forall \, p \in [p_1,p_2] \qquad \psi'(p) = (p-p_1)^{\rho_1 -1} (p_2-p)^{\rho_2 -1} \tilde{\psi}(p) \; ,
	\end{equation*}
	where $\tilde{\psi} \in \mathcal{C}^N\big([p_1,p_2], \R\big)$ is assumed positive, with $N \in \N \backslash \{0\}$. The points $p_j$ $(j=1,2)$ are called \emph{stationary points} of $\psi$ of order $\rho_j -1$.\\ \\
\textbf{Assumption (A$_{\mu_1,\mu_2,N}$).} For $0 < \mu_1 , \mu_2 \leqslant 1$, let $U : (p_1, p_2) \longrightarrow \C$ be a function defined by
	\begin{equation*}
		\forall \, p \in (p_1,p_2) \qquad U(p )= (p - p_1)^{\mu_1 -1} (p_2 - p)^{\mu_2 -1} \tilde{u}(p) \; ,
	\end{equation*}
	where $\tilde{u} \in \mathcal{C}^N\big([p_1,p_2], \C\big)$, with $N \in \N \backslash \{0\}$, and $\tilde{u}(p_j) \neq 0$ if $\mu_j \neq 1$ $(j=1,2)$. The points $p_j$ are called \emph{singularities} of $U$.\\

Now we define some objects that will be used throughout this paper.

\begin{def1} \label{DEF1}
	\begin{enumerate}
	\item Fix $q \in (p_1,p_2)$. For $j=1,2$, let $\varphi_j : I_j \longrightarrow \R$ be the functions defined by
	\begin{equation*}
		\varphi_1(p) := \big( \psi(p) - \psi(p_1) \big)^{1/\rho_1} \qquad \text{and} \qquad \varphi_2(p) := \big( \psi(p_2) - \psi(p) \big)^{1/\rho_2} \; ,
	\end{equation*}
	with $I_1 := [p_1, q]$, $I_2 := [q, p_2]$ and $s_1 := \varphi_1(q)$, $s_2 := \varphi_2(q)$.
	\item For $j=1,2$, define $k_j : (0,s_j] \longrightarrow \C$ by
	\begin{equation*}
			k_j(s) := U\big( \varphi_j^{-1}(s) \big) \, s^{1- \mu_j} \, \big( \varphi_j^{-1} \big)'(s) \; ,
	\end{equation*}
	which can be extended to the interval $[0,s_j]$ (see Proposition \ref{PROP4}).
	\item For $s > 0$, the complex curve $\Lambda^{(j)}(s)$ is defined by
	\begin{equation*}
		\Lambda^{(j)}(s) := \left\{ s + t e^{(-1)^{j+1} i \frac{\pi}{2\rho_j}} \, \Big| \, t \geqslant 0 \right\} \; .
	\end{equation*}
	\end{enumerate}
\end{def1}

\vspace{0.4cm}
	
\begin{thm1} \label{THM1}
	Let $N \in \N \backslash \{0\}$, let $\rho_1, \rho_2 \geqslant 1$ and $0 < \mu_1, \mu_2 < 1$. Suppose that the functions $\psi : [p_1,p_2] \longrightarrow \R$ and $U : (p_1, p_2) \longrightarrow \C$ satisfy Assumption \emph{(P$_{\rho_1,\rho_2,N}$)} and Assumption \emph{(A$_{\mu_1,\mu_2,N}$)}, respectively. Then for $j=1,2$, there are functions $A^{(j)}$, $R_1^{(j)}$, $R_2^{(j)} : (0,+\infty) \longrightarrow \C$ such that:
	\begin{equation*}
		\left\{ \begin{array}{rl}
			& \displaystyle \int_{p_1}^{p_2} U(p) e^{i \omega \psi(p)} \, dp = \sum_{j=1,2} A^{(j)}(\omega) + R_1^{(j)}(\omega,q) + R_2^{(j)}(\omega,q) \; , \\
			& \displaystyle \left| R_1^{(j)}(\omega,q) \right| \leqslant \frac{1}{\rho_j} \, \Gamma\bigg(\frac{1}{\rho_j}\bigg) \int_0^{s_j} s^{\mu_j-1} \big|(k_j)'(s) \big| \, ds \; \omega^{-\frac{1}{\rho_j}} \; , \\
			& \displaystyle \left| R_2^{(j)}(\omega,q) \right| \leqslant \frac{\rho_j-\mu_j}{\rho_j} \, \Gamma\bigg(\frac{1}{\rho_j}\bigg) \Big| U(q) \, (\varphi_j)'(q)^{-1} \Big| \varphi_j(q)^{-\rho_j} \, \omega^{-\left(1+\frac{1}{\rho_j}\right)} \; ,
		\end{array} \right.
	\end{equation*}
	for all $\omega > 0$ and for a fixed $q \in (p_1,p_2)$. For $j=1,2$ and $\omega > 0$, we have defined
	\begin{equation*}
		\begin{aligned}
			\bullet \quad & A^{(j)}(\omega) :=  e^{i \omega \psi(p_j)} \, k_j(0) \, \Theta^{(j)}(\rho_j,\mu_j) \, \omega^{-\frac{\mu_j}{\rho_j}} \; , \\
			\bullet \quad & R_1^{(j)}(\omega,q) := (-1)^{j} \, e^{i\omega \psi(p_j)} \int_0^{s_j} \phi^{(j)}(s,\omega, \rho_j, \mu_j) \, (k_j)'(s) \, ds \; , \\
			\bullet \quad & R_2^{(j)}(\omega,q) := (-1)^j \, i \, \frac{\mu_j-\rho_j}{\rho_j} \, e^{i \omega \psi(p_j)} \, k_j(s_j) \int_{\Lambda^{(j)}(s_j)} z^{\mu_j-\rho_j-1} e^{(-1)^{j+1} i \omega z^{\rho_j}} \, dz \; \omega^{-1} \; ,
		\end{aligned}
	\end{equation*}
	where
	\begin{equation*}
		\begin{aligned}
			\bullet \quad & \Theta^{(j)}(\rho_j,\mu_j) := \frac{(-1)^{j+1}}{\rho_j} \, \Gamma \left(\frac{\mu_j}{\rho_j}\right) \, e^{(-1)^{j+1} i \frac{\pi}{2} \, \frac{\mu_j}{\rho_j}} \; ,  \\
			\bullet \quad & \phi^{(j)}(s,\omega,\rho_j, \mu_j) := - \int_{\Lambda^{(j)}(s)} z^{\mu_j -1} e^{(-1)^{j+1} i\omega z^{\rho_j}} \, dz \; .
		\end{aligned}
	\end{equation*}
\end{thm1}

\begin{proof}
	For fixed $\rho_j \geqslant 1$ and $0 < \mu_j < 1$, we shall note $\phi^{(j)}(s,\omega)$ instead of $\phi^{(j)}(s,\omega,\rho_j,\mu_j)$. Now we choose $\omega > 0$ and we divide the proof in five steps.\\ \\
	\textit{First step: Splitting of the integral.} Fix a point $q \in (p_1,p_2)$ and split the integral at this point:
	\begin{equation*}
		\int_{p_1}^{p_2} U(p) e^{i \omega \psi(p)} \, dp = \tilde{I}^{(1)}(\omega,q) + \tilde{I}^{(2)}(\omega,q) \; ,
	\end{equation*}
	where
	\begin{equation*}
		\tilde{I}^{(1)}(\omega,q) := \int_{p_1}^{q} U(p) e^{i \omega \psi(p)} \, dp \quad , \quad \tilde{I}^{(2)}(\omega,q) := \int_{q}^{p_2} U(p) e^{i \omega \psi(p)} \, dp \; . \\ \\
	\end{equation*}
	\textit{Second step: Substitution.} Since Proposition \ref{PROP3} shows  that $\varphi_1 : [p_1,q] \longrightarrow [0,s_1]$ is a $\mathcal{C}^{N+1}$-diffeomorphism, we obtain by setting $s = \varphi_1(p)$,
	\begin{equation*}
		\begin{aligned}
			\tilde{I}^{(1)}(\omega,q)	& = \int_{p_1}^{q} U(p) \, e^{i \omega \psi(p)} \, dp \\
								& = e^{i \omega \psi(p_1)} \int_0^{s_1} U \big( \varphi_1^{-1}(s) \big) \, e^{i \omega s^{\rho_1}} (\varphi_1^{-1})'(s) \, ds \\
								& = e^{i \omega \psi(p_1)} \int_0^{s_1} U \big( \varphi_1^{-1}(s) \big) s^{1-\mu_1} (\varphi_1^{-1})'(s) \, s^{\mu_1-1} e^{i \omega s^{\rho_1}} \, ds \\
								& = e^{i \omega \psi(p_1)} \int_0^{s_1} k_1(s) \, s^{\mu_1-1} e^{i \omega s^{\rho_1}} \, ds \; ,
		\end{aligned}
	\end{equation*}
	where $k_1$ is given in Definition \ref{DEF1}. In a similar way, we obtain
	\begin{equation*}
		\tilde{I}^{(2)}(\omega,q) = - e^{i \omega \psi(p_2)} \int_0^{s_2} k_2(s) \, s^{\mu_2-1} e^{-i \omega s^{\rho_2}} \, ds \; ,
	\end{equation*}
	with $k_2$ defined above. Note that the minus sign comes from the decrease of $\varphi_2$.\\ \\
	\textit{Third step: Integration by parts.} Proposition \ref{COR1} furnishes a primitive of the function $s \longmapsto s^{\mu_j - 1} e^{(-1)^{j+1} i \omega s^{\rho_j}}$, and Proposition \ref{PROP4} assures that $k_j \in \mathcal{C}^1\big([0,s_1]\big)$. Thus we can integrate by parts:
	\begin{equation*}
		\begin{aligned}
			\tilde{I}^{(1)}(\omega,q)	& = e^{i \omega \psi(p_1)} \int_0^{s_1} k_1(s) \, s^{\mu_1-1} e^{i \omega s^{\rho_1}} \, ds \\
														& = \phi^{(1)}(s_1,\omega) \, k_1(s_1) \, e^{i \omega \psi(p_1)} - \phi^{(1)}(0,\omega) \, k_1(0) \, e^{i \omega \psi(p_1)}\\
							& \qquad \qquad - e^{i \omega \psi(p_1)} \int_0^{s_1} \phi^{(1)}(s,\omega) (k_1)'(s) \, ds \; .
		\end{aligned}
	\end{equation*}
	By a similar calculation, we have
	\begin{equation} \label{I_2}
		\begin{aligned}
			\tilde{I}^{(2)}(\omega,q)	& = \phi^{(2)}(0,\omega) \, k_2(0) \, e^{i \omega \psi(p_2)} - \phi^{(2)}(s_2,\omega) \, k_2(s_2) \, e^{i \omega \psi(p_2)}\\
								& \qquad \qquad + e^{i \omega \psi(p_2)} \int_0^{s_2} \phi^{(2)}(s,\omega) (k_2)'(s) \, ds \; .
		\end{aligned}
	\end{equation}
	
	\noindent \textit{Fourth step: Cancellation.} The aim of this step is to simplify the difference:
	\begin{equation*}
		\phi^{(1)}(s_1,\omega) \, k_1 (s_1) \, e^{i \omega \psi(p_1)} - \phi^{(2)}(s_2,\omega) \, k_2 (s_2) \, e^{i \omega \psi(p_2)} \; .
	\end{equation*}
	The functions $s \mapsto \phi^{(j)}(s,\omega)$ are given by oscillatory integrals. We can expand them with respect to $\omega$ and we show that the first terms cancel each other out.\\ Since $s_j > 0$, we note that the derivative of the function $t \longmapsto \big(s_j+t e^{(-1)^{j+1} i \frac{\pi}{2\rho_j}}\big)^{\rho_j}$ does not vanish for all $t > 0$  and one can write:
	\begin{equation*}
		e^{i \omega \left(s_1+t e^{i \frac{\pi}{2\rho_1}}\right)^{\rho_1}} = e^{-i \frac{\pi}{2\rho_1}} \, (i \omega \rho_1)^{-1} \left(s_1+t e^{i \frac{\pi}{2\rho_1}}\right)^{1-\rho_1} \frac{d}{dt} \left[e^{i \omega \left(s_1+t e^{i \frac{\pi}{2\rho_1}}\right)^{\rho_1}}\right] \; .
	\end{equation*}
	Putting this equality for $j=1$ in the definition of $\phi^{(1)}(s_1,\omega)$ and carrying out an integration by parts lead to
	\begin{equation} \label{expphi1}
		\begin{aligned}
			\phi^{(1)}(s_1,\omega)	& = - (i \omega \rho_1)^{-1} \int_0^{+\infty} \left(s_1+te^{i \frac{\pi}{2 \rho_1}}\right)^{\mu_1-\rho_1} \frac{d}{dt}\left[e^{i \omega \left(s_1+te^{i \frac{\pi}{2 \rho_1}}\right)^{\rho_1}}\right] \, dt \\
														& = (i \omega \rho_1)^{-1} (s_1)^{\mu_1-\rho_1} \, e^{i \omega (s_1)^{\rho_1}} \\
														& \qquad + \frac{\mu_1-\rho_1}{i \omega \rho_1} \, e^{i \frac{\pi}{2 \rho_1}} \int_0^{+\infty} \left(s_1+te^{i \frac{\pi}{2 \rho_1}}\right)^{\mu_1-\rho_1-1} e^{i \omega \left(s_1+te^{i \frac{\pi}{2 \rho_1}}\right)^{\rho_1}} \, dt \; .
		\end{aligned}
	\end{equation}
	We remark that the boundary term at infinity is $0$; indeed, we observe that
	\begin{equation*}
		s_1 \leqslant \left| s_1+te^{i \frac{\pi}{2 \rho_1}} \right| \qquad \Longrightarrow \qquad (s_1)^{\mu_1-\rho_1} \geqslant \left| s_1+te^{i \frac{\pi}{2 \rho_1}} \right|^{\mu_1-\rho_1} \; ,
	\end{equation*}
	because $\mu_1 < \rho_1$, and by using Proposition \ref{PROP2}, we obtain
	\begin{equation*}
		\forall \, t > 0 \qquad \left| \left(s_1+te^{i \frac{\pi}{2 \rho_1}}\right)^{\mu_1-\rho_1} e^{i \omega \left(s_1+te^{i \frac{\pi}{2 \rho_1}}\right)^{\rho_1}} \right| \leqslant (s_1)^{\mu_1-\rho_1} e^{-\omega t^{\rho_1}} \; \longrightarrow \; 0 \quad , \quad t \rightarrow + \infty \; .
	\end{equation*}
	In a similar way, we obtain
	\begin{equation} \label{expphi2}
		\begin{aligned}
			\phi^{(2)}(s_2,\omega)	& = -(i \omega \rho_2)^{-1} (s_2)^{\mu_2-\rho_2} \, e^{-i \omega (s_2)^{\rho_2}} \\
														& \qquad - \frac{\mu_2-\rho_2}{i \omega \rho_2} \, e^{-i \frac{\pi}{2 \rho_2}} \int_0^{+\infty} \left(s_2+te^{-i \frac{\pi}{2 \rho_2}}\right)^{\mu_2-\rho_2-1} e^{-i \omega \left(s_2+te^{-i \frac{\pi}{2 \rho_2}}\right)^{\rho_2}} \, dt \; .
		\end{aligned}
	\end{equation}
	Furthermore, by the definitions of $k_j$ and $s_j := \varphi_j(q)$, we get
	\begin{equation*}
		k_j(s_j) = U\big( \varphi_j^{-1}(s_j) \big) \, s_j^{1-\mu_j} \, (\varphi_j^{-1})'(s_j) = U(q) \, \varphi_j(q)^{1-\mu_j} \, (\varphi_j)'(q)^{-1} \; .
	\end{equation*}
	Now we multiply the last expression by the expansion of $\phi^{(1)}(s_1,\omega)$ found in \eqref{expphi1},
	\begin{equation} \label{simpphi1}
		\begin{aligned}
			\phi^{(1)} \, (s_1,\omega) & \, k_1(s_1) e^{i \omega \psi(p_1)} = (i \omega \rho_1)^{-1} e^{i \omega \left((s_1)^{\rho_1} + \psi(p_1)\right)} \, U(q) \, \varphi_1(q)^{1-\rho_1} \, (\varphi_1)'(q)^{-1} \\
									& - i \, \frac{\mu_1 - \rho_1}{\rho_1} \, e^{i \omega \psi(p_1)} \, e^{i \frac{\pi}{2 \rho_1}} \, k_1(s_1) \int_0^{+\infty} \left(s_1+te^{i \frac{\pi}{2 \rho_1}}\right)^{\mu_1-\rho_1-1} e^{i \omega \left(s_1+te^{i \frac{\pi}{2 \rho_1}}\right)^{\rho_1}} \, dt \; \omega^{-1} \; .
		\end{aligned}
	\end{equation}
	The definition of $s_1$ gives $(s_1)^{\rho_1} = \psi(q) - \psi(p_1)$ and by the regularity of $\varphi_1$, one has
	\begin{equation*}
		\rho_1 \, (\varphi_1)'(q) \, \varphi_1(q)^{\rho_1-1} = \frac{d}{dp}\Big[ (\varphi_1)^{\rho_1} \Big](q) = \psi'(q) \; ;
	\end{equation*}
	hence the first term can be simplified; moreover the integral in \eqref{simpphi1} can be written as an integral on the curve $\Lambda^{(1)}(s_1)$ in the complex plane. These considerations lead to
	\begin{equation} \label{phi1}
		\begin{aligned}
			\phi^{(1)}(s_1,\omega) \, & k_1(s_1)  \, e^{i \omega \psi(p_1)} = -i \omega^{-1} \, e^{i \omega \psi(q)} \, \frac{U(q)}{\psi'(q)} \\
										& \qquad - i \, \frac{\mu_1 - \rho_1}{\rho_1}  \, e^{i \omega \psi(p_1)} \, k_1(s_1)  \int_{\Lambda^{(1)}(s_1)} z^{\mu_1-\rho_1-1} e^{i \omega z^{\rho_1}} \, dz \; \omega^{-1} \; .
		\end{aligned}
	\end{equation}
	In a similar way, we obtain
	\begin{equation} \label{phi2}
		\begin{aligned}
			\phi^{(2)} \, (s_2,\omega) & \, k_2(s_2) e^{i \omega \psi(p_2)}	= -i \omega^{-1} \, e^{i \omega \psi(q)} \, \frac{U(q)}{\psi'(q)} \\
										& \qquad + i \, \frac{\mu_2 - \rho_2}{\rho_2} \, e^{i \omega \psi(p_2)}  k_2(s_2) \int_{\Lambda^{(2)}(s_2)} z^{\mu_2-\rho_2-1} e^{-i \omega z^{\rho_2}} \, dz \; \omega^{-1} \; .
		\end{aligned}
	\end{equation}
	Now remark that if we add $\tilde{I}^{(1)}(\omega,q)$ and $\tilde{I}^{(2)}(\omega,q)$, then the two first terms of \eqref{phi1} and \eqref{phi2} cancel each other out. Consequently, we are able to write the initial integral as follows:
	\begin{equation*}
		\begin{aligned}
			\int_{p_1}^{p_2} U(p) e^{i \omega \psi(p)} dp	& = - \phi^{(1)}(0,\omega) \, k_1(0) \, e^{i \omega \psi(p_1)} - e^{i \omega \psi(p_1)} \int_0^{s_1} \phi^{(1)}(s,\omega) (k_1)'(s) \, ds \\
							& \qquad \qquad - i \, \frac{\mu_1-\rho_1}{\rho_1}  \, e^{i \omega \psi(p_1)} \, k_1(s_1)  \int_{\Lambda^{(1)}(s_1)} z^{\mu_1-\rho_1-1} e^{i \omega z^{\rho_1}} \, dz \; \omega^{-1} \\
							& \quad + \phi^{(2)}(0,\omega) \, k_2(0) \, e^{i \omega \psi(p_2)} + e^{i \omega \psi(p_2)} \int_0^{s_2} \phi^{(2)}(s,\omega) (k_2)'(s) \, ds \\
							& \qquad \qquad + i \, \frac{\mu_2-\rho_2}{\rho_2} \, e^{i \omega \psi(p_2)} \, k_2(s_2) \int_{\Lambda^{(2)}(s_2)} z^{\mu_2-\rho_2-1} e^{-i \omega z^{\rho_2}} \, dz \; \omega^{-1} \\
							& =: \sum_{j=1,2} A^{(j)}(\omega) + R_1^{(j)}(\omega,q) + R_2^{(j)}(\omega,q) \; .
		\end{aligned}
	\end{equation*}
	According to \cite[Theorem 1.3]{article1}, we have $\phi^{(j)}(0,\omega) = \Theta^{(j)}(\rho_j,\mu_j) \, \omega^{-\frac{\mu_j}{\rho_j}}$, where we put
	\begin{equation*}
		\Theta^{(j)}(\rho_j,\mu_j) := \frac{(-1)^{j+1}}{\rho_j} \, \Gamma \left(\frac{\mu_j}{\rho_j}\right) \, e^{(-1)^{j+1} i \frac{\pi}{2} \, \frac{\mu_j}{\rho_j}} \; .
	\end{equation*}
	Thus we define $A^{(j)}(\omega)$ as follows:
	\begin{equation*}
		A^{(j)}(\omega) := (-1)^{j} \phi^{(j)}(0,\omega) \, k_j(0) \, e^{i \omega \psi(p_j)} = 	e^{i \omega \psi(p_j)} \, k_j(0) \, \Theta^{(j)}(\rho_j,\mu_j) \, \omega^{-\frac{\mu_j}{\rho_j}} \; . \\
	\end{equation*}
	
	\noindent \textit{Fifth step: Remainder estimates.} To conclude the proof, we have to establish estimates for the remainder terms. Using \cite[Theorem 1.3]{article1}, we obtain
	\begin{equation*}
		\left| e^{i \omega \psi(p_j)} \int_0^{s_j} \phi^{(j)}(s,\omega, \rho_j, \mu_j) \, (k_j)'(s) \, ds \right| \leqslant \frac{1}{\rho_j} \, \Gamma\left(\frac{1}{\rho_j}\right) \int_0^{s_j} s^{\mu_j-1} \big| (k_1)'(s) \big| ds \, \omega^{-\frac{1}{\rho_j}} \; .
	\end{equation*}
	Note that we do not have to deal with any cut-off function in this situation. Now let us estimate $R_2^{(j)}(\omega,q)$:
	\begin{align}
		\Bigg| i \, \frac{\mu_j-\rho_j}{\omega \rho_j}	& \, e^{i \omega \psi(p_j)} \, k_j(s_j) \, e^{(-1)^{j+1} i \frac{\pi}{2 \rho_j}} \nonumber \\
												& \qquad \times \int_0^{+\infty} \left(s_j+te^{(-1)^{j+1} i \frac{\pi}{2\rho_j}}\right)^{\mu_j-\rho_j-1} e^{(-1)^{j+1} i \omega \left(s_j+te^{(-1)^{j+1}i \frac{\pi}{2 \rho_j}}\right)^{\rho_j}} \, dt \, \Bigg| \nonumber \\
												& \label{ineq1} \leqslant \frac{\rho_j-\mu_j}{\rho_j} \, \big| k_j(s_j) \big| \, \omega^{-1} \,  (s_j)^{\mu_j-\rho_j-1} \int_0^{+\infty} \left| e^{(-1)^{j+1} i \omega \left(s_j+te^{(-1)^{j+1}i \frac{\pi}{2 \rho_j}}\right)^{\rho_j}} \right| \, dt \\
												& \label{ineq2} \leqslant \frac{\rho_j-\mu_j}{\rho_j} \, \Big| U \big( \varphi_j^{-1}(s_j)\big) \, (\varphi_j^{-1})'(s_j) \Big| \, (s_j)^{-\rho_j} \, \omega^{-1} \int_0^{+\infty} e^{-\omega t^{\rho_j}} \, dt \\
												& \label{ineq3} = \frac{\rho_j-\mu_j}{\rho_j} \, \Gamma\left(\frac{1}{\rho_j}\right) \, \big| U(q) \, (\varphi_j)'(q)^{-1} \big| \, \varphi_j(q)^{-\rho_j} \,  \omega^{-\left(1+\frac{1}{\rho_j}\right)} \; ;
	\end{align}
	\begin{itemize}
		\item \eqref{ineq1}: $\displaystyle s_j \leqslant \left| s_j+te^{(-1)^{j+1}i \frac{\pi}{2 \rho_j}} \right|$ ;
		\item \eqref{ineq2}: definition of the function $k_j$ and Proposition \ref{PROP2} ;
		\item \eqref{ineq3}: $\displaystyle \int_0^{+\infty} e^{-\omega t^{\rho_j}} dt = \Gamma\left(\frac{1}{\rho_j} \right) \omega^{-\frac{1}{\rho_j}}$ and $q = \varphi_j^{-1}(s_j)$.
	\end{itemize}
	We remark finally that the decay rates of $A^{(j)}(\omega)$, $R_1^{(j)}(\omega,q)$ and $R_2^{(j)}(\omega,q)$ are $\omega^{-\frac{\mu_j}{\rho_j}}$, $\omega^{-\frac{1}{\rho_j}}$ and $\omega^{-\big(1+\frac{1}{\rho_j}\big)}$, respectively. Thus the decay rates of the remainder related to $p_j$ and the remainder related to $q$ are higher to the one of the first term related to $p_j$. This ends the proof.
\end{proof}

\vspace{0.4cm}

For fixed $q \in (p_1,p_2)$, we notice that if $\mu_j=1$ (for $j=1,2$), then the decay rate of $R_1^{(j)}(\omega,q)$ and $A^{(j)}(\omega)$ with respect to $\omega$ are the same. So we shall use the estimate provided by \cite[Theorem 1.7]{article1} to obtain a better decay rate. Furthermore we observe that $R_2^{(j)}(\omega,q)$ is always negligible as compared with $A^{(j)}(\omega)$ even if $\mu_j=1$.

\begin{thm2} \label{THM2}
	Let $N \in \N \backslash \{0\}$, $\rho_1, \rho_2 \geqslant 2$ and assume  $\mu_j = 1$ for a certain $j \in \{1,2\}$. Suppose that the functions $\psi : [p_1,p_2] \longrightarrow \R$ and $U : (p_1, p_2) \longrightarrow \C$ satisfy Assumption \emph{(P$_{\rho_1,\rho_2,N}$)} and Assumption \emph{(A$_{\mu_1,\mu_2,N}$)}, respectively. Then the statement of Theorem \ref{THM1} is still true and, for $\gamma \in (0,1)$ and $\delta := \rho_j^{-1}(\gamma + 1) \in \big( \frac{1}{\rho_j}, \frac{2}{\rho_j} \big)$, we have more precise estimates for the remainder related to $p_j$, namely:
	\begin{equation} \label{rest2}
		\left| R_1^{(j)}(\omega,q) \right| \leqslant L_{\gamma, \rho_j} \int_0^{s_j} s^{-\gamma} \big| (k_j)'(s) \big| \, ds \; \omega^{-\delta} \; ,
	\end{equation}
	for all $\omega > 0$. See \cite[Lemmas 1.4, 1.6]{article1} for the definition of the constant $L_{\gamma, \rho_j} > 0$.
\end{thm2}

\begin{proof}
	We obtain the above estimates following the lines of the proof of \cite[Theorem 1.7]{article1}.
\end{proof}


\section{Approaching stationary points and amplitude singularities:
the first term and error between blow-up and decay}

\hspace{2.5ex} In this section, we consider a family of oscillatory integrals with respect to a large parameter $\omega$. We suppose that the phase function is a polynomial of degree $2$ and has its stationary point $p_0$ inside $(p_1,p_2)$, which contains the support of the amplitude. We suppose in addition that the amplitude has a singularity at $p_1$, the left endpoint of the  integration interval. The aim of this section is to furnish explicit estimates of the difference of the oscillatory integrals and the first term of their expansions.

In our first result, we establish this estimate in the whole space of the parameters $\omega$ and $p_0$. We find that the first term of the expansion, as well as the error estimate, blows up when the stationary point $p_0$ approaches the singularity $p_1$.


In the second result, we suppose that $p_0$ approaches $p_1$ with a certain convergence speed, described by the parameter $\varepsilon > 0$, as the large parameter $\omega$ tends to infinity. Below a certain threshold, the convergence speed is sufficiently slow so that the decay with respect to $\omega$ compensates the blow-up. Using the preciseness of the preceding result, we furnish finally asymptotic expansions of the oscillatory integrals on curves in the space of the parameters.\\

\begin{thm3} \label{THM3}
	Let $p_1,p_2$ be two real numbers such that $-\infty < p_1 < p_2 < +\infty$. Let $p_0, c \in \R$ be two parameters and define $\psi : [p_1,p_2] \longrightarrow \R$ by
	\begin{equation*}
		\psi(p) := - (p-p_0 )^2 + c \; .
	\end{equation*}
	Suppose that $p_0 \in (p_1, p_2)$ and define the following integrals for all $\omega > 0$:
	\begin{equation*}
		I^{(1)}(\omega,p_0) := \int_{p_1}^{p_0} U(p) \, e^{i \omega \psi(p)} \, dp \qquad , \qquad I^{(2)}(\omega,p_0) := \int_{p_0}^{p_2} U(p) \, e^{i \omega \psi(p)} \, dp \; ,
	\end{equation*}
	where $U$ satisfies Assumption \emph{(A$_{\mu,1,1}$)} on $[p_1,p_2]$ with $\mu \in (0,1)$ and $U(p_2)=0$. Then
	\begin{itemize}
		\item there exist $\tilde{K}_{\mu}^{(1)}(\omega,p_0,U), \tilde{H}^{(1)}(\omega,p_0,U) \in \C$ and constants $R_k^{(1)}(U) \geqslant 0$ $(k=1,...,6)$ independent on $\omega$ and $p_0$, such that the following estimate holds:
		\begin{equation*}
			\begin{aligned}
				\Big| I^{(1)}(\omega,p_0) - \tilde{K}_{\mu}^{(1)}(\omega,p_0,U) \, (p_0 - p_1)^{-\mu} \, \omega^{-\mu} -	& \tilde{H}^{(1)}(\omega,p_0,U) \, (p_0 - p_1)^{\mu - 1} \, \omega^{-\frac{1}{2}} \Big| \\
											& \leqslant \sum_{k = 1}^6 \, R_k^{(1)}(U) \, (p_0 - p_1)^{-\alpha_k^1} \, \omega^{-\alpha_k^2} \; ,
			\end{aligned}
		\end{equation*}
		where $\alpha_k^1, \alpha_k^2 \geqslant 0$;
		\item there exist $\tilde{H}^{(2)}(\omega,p_0,U) \in \C$ and constants $R_k^{(2)}(U) \geqslant 0$ $(k=1,2)$ independent on $\omega$ and $p_0$, such that the following estimate holds:
		\begin{equation*}
			\Big| I^{(2)}(\omega,p_0) -	 \tilde{H}^{(2)}(\omega,p_0,U) \, (p_0 - p_1)^{\mu - 1} \, \omega^{-\frac{1}{2}} \Big| \leqslant \sum_{k = 1}^2 \, R_k^{(2)}(U) \, (p_0 - p_1)^{-\beta_k^1} \, \omega^{-\beta_k^2} \; ,
		\end{equation*}
		where $\beta_k^1, \beta_k^2 \geqslant 0$.
	\end{itemize}
\end{thm3}

\begin{rem2}
	\em
	\begin{enumerate}
		\item If we suppose $p_0 - p_1 > \varepsilon$, for a certain $\varepsilon > 0$, then we obtain uniform estimates in $p_0$ of the remainder. Applying this to the free Schrödinger equation leads to the space-time cones introduced in \cite{article1}.
		\item The preciseness that we obtain in this paper shows that the remainder estimates does not blow up when the stationary point $p_0$ tends to the regular point $p_2$. In \cite{article1}, the smooth cut-off function caused an artificial blow-up as $p_0$ tends to $p_2$, which was not controllable.
		\item A certain number $\alpha_k^1$ can be chosen negative. Without loss of generality, we can suppose that it is positive for simplicity.
	\end{enumerate}
\end{rem2}

\begin{proof}
	\noindent \textit{Study of $I^{(1)}(\omega,p_0)$}. For all $p \in [p_1,p_0]$, we have:
	\begin{equation*}
		\psi'(p) = 2 (p_0-p) \; .
	\end{equation*}
	By setting $\tilde{\psi} = 2$, we observe that $\psi$ verifies Assumption \emph{(P$_{1,2,N}$)} on $[p_1,p_0]$, for all $N \geqslant 1$. Then we can apply Theorem \ref{THM1} to $I^{(1)}(\omega,p_0)$; here we choose arbitrarily $q := q(p_0) = p_1 + \frac{p_0-p_1}{2} = p_0 - \frac{p_0-p_1}{2}$ for simplicity. Then we get
	\begin{equation*}
		I^{(1)}(\omega,p_0) = \sum_{j=1,2} A^{(j)}(\omega,p_0) + R_1^{(j)}(\omega,p_0) + R_2^{(j)}(\omega,p_0) \; ,
	\end{equation*}
	with
	\begin{equation*}
		\begin{aligned}
			& \bullet \quad A^{(1)}(\omega,p_0) = e^{i \omega \psi(p_1)} \, k_1(0) \, \Theta^{(1)}(1,\mu) \, \omega^{-\mu} = \Gamma(\mu) \, e^{i \frac{\pi \mu}{2}} \, e^{i \omega \psi(p_1)} \, k_1(0) \, \omega^{-\mu} \; , \\
			& \bullet \quad A^{(2)}(\omega,p_0) = e^{i \omega \psi(p_0)} \, k_2(0) \, \Theta^{(2)}(2,1) \, \omega^{-\frac{1}{2}} = - \frac{\sqrt{\pi}}{2} \, e^{-i \frac{\pi}{4}} \, e^{i \omega \psi(p_0)} \, k_2(0) \, \omega^{-\frac{1}{2}} \; .
		\end{aligned}
	\end{equation*}
	To clarify these terms, let us study the functions $(\varphi_1^{-1})'$ and $(\varphi_2^{-1})'$. On the one hand, by the simple definition of $\varphi_1$, we obtain
	\begin{equation} \label{phi122}
		\varphi_1(p) = \psi(p) - \psi(p_1) \qquad \Longrightarrow \qquad (\varphi_1)'(p) = \psi'(p) = 2 (p_0-p) \; ,
	\end{equation}
	for all $p \in [p_1,q]$. On the other hand, by the definition of $\varphi_2$ and the expression of $\psi$, one has
	\begin{equation*}
		\varphi_2(p) = \big( \psi(p_0) - \psi(p) \big)^{\frac{1}{2}} = (p_0-p) \; ,
	\end{equation*}
	for every $p \in [q,p_0]$. So $(\varphi_2^{-1})'(s) = -1$ and $(\varphi_2^{-1})''(s) = 0$ for all $s \in [0,s_2]$. Then we compute $k_1(0)$ by using the following representation of $k_1$ (see \cite[(35)]{article1}):
	\begin{align}
		\label{k_1} k_1(s)	& = \left( \int_0^1 (\varphi_1^{-1})'(sy) \, dy \right)^{\mu-1} \tilde{u}\big(\varphi_1^{-1}(s) \big) (\varphi_1^{-1})'(s) \\
				& \underset{s \longrightarrow 0^+}{\longrightarrow} \; (\varphi_1^{-1})'(0)^{\mu} \, \tilde{u}\big( \varphi_1^{-1}(0) \big) = \psi'(p_1)^{-\mu} \, \tilde{u}(p_1) = \big( 2 (p_0 - p_1) \big)^{-\mu} \, \tilde{u}(p_1) \; . \nonumber
	\end{align}
	To compute $k_2(0)$, we use the definition of $k_2$ because $p_0$ is not a singularity of the amplitude:
	\begin{equation*}
		k_2(s) = U\big( \varphi_2^{-1}(s) \big) (\varphi_2^{-1})'(s) = -U\big( \varphi_2^{-1}(s) \big) \; \underset{s \longrightarrow 0^+}{\longrightarrow} \; -U(p_0) = -(p_0 - p_1)^{\mu-1} \, \tilde{u}(p_0) \; .
	\end{equation*}
	Therefore we obtain
	\begin{align}
		& \label{B11} \bullet \quad A^{(1)}(\omega,p_0) = \frac{\Gamma(\mu)}{2^{\mu}} \, e^{i \frac{\pi \mu}{2}} \, e^{i \omega \psi(p_1)} \, \tilde{u}(p_1) \, (p_0 - p_1)^{-\mu} \, \omega^{-\mu} \; , \\
		& \label{B12} \bullet \quad A^{(2)}(\omega,p_0) = \frac{\sqrt{\pi}}{2} \, e^{-i \frac{\pi}{4}} \, e^{i \omega \psi(p_0)} \, \tilde{u}(p_0) \, (p_0 - p_1)^{\mu-1} \, \omega^{-\frac{1}{2}} \; .
	\end{align}
	
	Now let us control precisely the remainder terms. To do so, we need to study the functions $k_j$ and $\varphi_j$. Firstly, employing the relations $\frac{1}{2}(p_0-p_1) = p_0-q = q-p_1$ and \eqref{phi122} leads to
	\begin{equation*}
		(p_0-p_1) = 2(p_0-q) \leqslant (\varphi_1)'(p) \leqslant 2 (p_0-p_1) \; ,
	\end{equation*}
	for all $p \in [p_1,q]$. It follows
	\begin{equation} \label{derinvphi1}
		\forall \,  s \in [0,s_1] \qquad \big(2 (p_0-p_1) \big)^{-1} \leqslant (\varphi_1^{-1})'(s) \leqslant (p_0 - p_1)^{-1} \; .
	\end{equation}
	Moreover by the equality $(\varphi_1^{-1})''(s) = - (\varphi_1)''\big(\varphi_1^{-1}(s)\big) \, (\varphi_1^{-1})'(s)^3$, we have
	\begin{equation} \label{phisec}
		\forall \, s \in [0,s_1] \qquad (\varphi_1^{-1})''(s) = 2 \, (\varphi_1^{-1})'(s)^3 \leqslant 2 \, (p_0-p_1)^{-3} \; .
	\end{equation}
	Furthermore by the mean value Theorem, it is possible to bound $s_1 = \varphi_1(q) = \varphi_1(q) - \varphi(p_1)$ from above and below:
	\begin{equation} \label{s1}
		\frac{1}{2} \, (p_0-p_1)^2 = (p_0-p_1)(q-p_1) \leqslant s_1 \leqslant 2(p_0-p_1) (q-p_1) = (p_0-p_1)^2 \; .
	\end{equation}
	Concerning $s_2$, we have the following expression:
	\begin{equation} \label{s2}
		s_2 = \varphi_2(q) = (p_0 - q) = \frac{1}{2} \, (p_0-p_1) \; .
	\end{equation}
	Now we study the functions $k_j$. For this purpose, we shall use the expression of $k_1$ given in \eqref{k_1}. Since $\psi$ satisfies Assumption (P$_{1,2,N}$) on $[p_1,p_0]$ for all $N \geqslant 1$, it follows that $\varphi_1$ is a $\mathcal{C}^{N+1}$-diffeomorphism by Proposition \ref{PROP3}. Thus one has the ability to differentiate under the integral sign the function $s \longmapsto \int_0^1 (\varphi_1^{-1})'(sy) dy$. Then the derivative of $k_1$ on $[0,s_1]$ is given by
	\begin{equation*}
		\begin{aligned}
			(k_1)'(s) =	& (\mu-1) \left( \int_0^1 y (\varphi_1^{-1})''(sy) \, dy \right) \left( \int_0^1 (\varphi_1^{-1})'(sy) \, dy \right)^{\mu-2} \tilde{u}\big(\varphi_1^{-1}(s) \big) (\varphi_1^{-1})'(s) \\
						& \qquad + \left( \int_0^1 (\varphi_1^{-1})'(sy) \, dy \right)^{\mu-1} \tilde{u}'\big(\varphi_1^{-1}(s) \big) (\varphi_1^{-1})'(s)^2 \\
						&  \qquad \qquad + \left( \int_0^1 (\varphi_1^{-1})'(sy) \, dy \right)^{\mu-1} \tilde{u}\big(\varphi_1^{-1}(s) \big) (\varphi_1^{-1})''(s) \; .
		\end{aligned}
	\end{equation*}
	In consequences, we obtain the following estimate:
	\begin{align}
		\big\| (k_1)' \big\|_{L^{\infty}(0,s_1)} & \leqslant \frac{1-\mu}{2} \, 2 \, (p_0-p_1)^{-3} \, \big(2 (p_0-p_1) \big)^{2-\mu} \, \| \tilde{u} \|_{L^{\infty}(p_1,p_2)} (p_0 - p_1)^{-1} \nonumber \\
												& \qquad + \big(2 (p_0-p_1) \big)^{1-\mu} \, \| \tilde{u}' \|_{L^{\infty}(p_1,p_2)} (p_0 - p_1)^{-2} \nonumber \\
												& \qquad \qquad + \big(2 (p_0-p_1) \big)^{1-\mu} \| \tilde{u} \|_{L^{\infty}(p_1,p_2)} \, 2 (p_0-p_1)^{-3} \nonumber \\
												& \label{derk1} \leqslant 2^{1-\mu} \, \| \tilde{u} \|_{W^{1,\infty}(p_1,p_2)} \, \Big( 2(2-\mu)(p_0-p_1)^{-(2+\mu)} + (p_0 - p_1)^{-(1+\mu)} \Big) \; .
	\end{align}
	
	As above, we study the function $k_2$ by employing its definition. Recalling that $U(p) = (p-p_1)^{\mu-1} \tilde{u}(p)$, we differentiate $k_2$:
	\begin{equation*}
		(k_2)'(s) = \Big( (\mu-1) \big( \varphi_2^{-1}(s) - p_1\big)^{\mu-2} \, \tilde{u}\big( \varphi_2^{-1}(s) \big) + \big( \varphi_2^{-1}(s) - p_1 \big)^{\mu-1} \, \tilde{u}'\big( \varphi_2^{-1}(s) \big) \Big) (\varphi_2^{-1})'(s)^2 \; ,
	\end{equation*}
	for $s \in [0,s_2]$. Use the fact that $\varphi_2^{-1}(s) \in [q,p_2]$ for $s \in [0,s_2]$ and the equality $(\varphi_2^{-1})'(s) = -1$ to obtain
	\begin{align}
		\big\| (k_2)' \big\|_{L^{\infty}(0,s_2)}	& \leqslant \Big((1-\mu) \, 2^{2-\mu} (p_0-p_1)^{\mu-2} \, \| \tilde{u} \|_{L^{\infty}(p_1,p_2)} + 2^{1-\mu} (p_0-p_1)^{\mu-1} \| \tilde{u}' \|_{L^{\infty}(p_1,p_2)} \Big) \nonumber \\
													& \label{k2} \leqslant 2^{1-\mu} \, \| \tilde{u} \|_{W^{1,\infty}(p_1,p_2)} \Big( 2(1-\mu) (p_0-p_1)^{\mu-2} + (p_0-p_1)^{\mu-1} \Big) \; .
	\end{align}
	These considerations permit to estimate the four remainders.
	\begin{itemize}
		\item \emph{Estimate of $R_1^{(1)}(\omega,p_0)$.} Theorem \ref{THM1} furnishes an estimate of $R^{(1)}(\omega,p_0)$. We combine it with the estimates of $(k_1)'$ \eqref{derk1} and $s_1$ \eqref{s1}:
		\begin{equation} \label{R111}
			\begin{aligned}
				\left| R_1^{(1)}(\omega,p_0) \right|	& \leqslant \int_0^{s_1} s^{\mu-1} \big|(k_1)'(s) \big| \, ds \; \omega^{-1} \\
														& \leqslant \mu^{-1} \, (s_1)^{\mu} \, \big\| (k_1)' \big\|_{L^{\infty}(0,s_1)} \, \omega^{-1} \\
														& \leqslant \frac{2^{1-\mu}}{\mu} \, \| \tilde{u} \|_{W^{1,\infty}(p_1,p_2)} \Big( 2(2-\mu)(p_0-p_1)^{\mu-2} + (p_0 - p_1)^{\mu-1} \Big) \omega^{-1} \; .
			\end{aligned}
		\end{equation}
		\item \emph{Estimate of $R_2^{(1)}(\omega,p_0)$.} The estimate of $R_2^{(1)}(\omega,p_0)$ from Theorem \ref{THM1} provides
		\begin{align}
			\left| R_2^{(1)}(\omega,p_0) \right|	& \leqslant (1-\mu) \, \left| U(q) \, (\varphi_1)'(q)^{-1} \right| \, \varphi_1(q)^{-1} \, \omega^{-2} \nonumber \\
													& \leqslant \frac{1-\mu}{2^{\mu-1}} \, \| \tilde{u} \|_{L^{\infty}(p_1,p_2)} (p_0 - p_1)^{\mu-1} \, (p_0-p_1)^{-1} \, \big( 2^{-1} (p_0-p_1)^2 \big)^{-1} \, \omega^{-2} \nonumber \\
													& \label{R211} = \frac{1-\mu}{2^{\mu-2}} \, \| \tilde{u} \|_{L^{\infty}(p_1,p_2)} (p_0 - p_1)^{\mu-4} \, \omega^{-2} \; ,
		\end{align}
		where the definition of $U$, inequalities \eqref{derinvphi1} and \eqref{s1} were used.
		\item \emph{Estimate of $R_1^{(2)}(\omega,p_0)$.} Here $\mu_2 = 1$, so we have to employ the estimate of $R_1^{(2)}(\omega,p_0)$ from Theorem \ref{THM2},
		\begin{equation} \label{R121}
			\begin{aligned}
				\left| R_1^{(2)}(\omega,p_0) \right|	& \leqslant L_{\gamma, 2} \int_0^{s_2} s^{-\gamma} \big| (k_2)'(s) \big| \, ds \; \omega^{-\delta} \\
														& \leqslant \frac{L_{\gamma, 2}}{1-\gamma} \, (s_2)^{1-\gamma} \big\| (k_2)' \big\|_{L^{\infty}(0,s_2)} \omega^{-\delta} \\
														& \leqslant \frac{L_{\gamma, 2}}{1-\gamma} \, 2^{\gamma-\mu} \, \| \tilde{u} \|_{W^{1,\infty}(p_1,p_2)} \Big( 2(1-\mu) (p_0-p_1)^{\mu-1-\gamma} + (p_0-p_1)^{\mu-\gamma} \Big) \omega^{-\delta} \; ,
			\end{aligned}
		\end{equation}
		the last inequality was obtained by using \eqref{s2} and \eqref{k2}. Since we want the exponents of $(p_0-p_1)$ to be negative, we require $\mu - \gamma \leqslant 0$; this is equivalent to $\frac{1}{2} < \frac{\mu + 1}{2} \leqslant \delta < 1$ by the relation $\gamma = 2 \delta - 1$.
		\item \emph{Estimate of $R_2^{(2)}(\omega,p_0)$.} We employ Theorem \ref{THM1} once again to control $R_2^{(2)}(\omega,p_0)$. The definition of $U$, the relation $(\varphi_2^{-1})' = -1$ and equality \eqref{s2} lead to
		\begin{align}
			\left| R_2^{(2)}(\omega,p_0) \right|	& \leqslant \frac{1}{2} \, \Gamma \left(\frac{1}{2}\right) \, \left| U(q) \, (\varphi_2)'(q)^{-1} \right| \, \varphi_2(q)^{-2} \, \omega^{-\frac{3}{2}} \nonumber \\
												& \leqslant \frac{\sqrt{\pi}}{2^{\mu}} \, \| \tilde{u} \|_{L^{\infty}(p_1,p_2)} \, (p_0 - p_1)^{\mu-1} \, \big(2^{-1} (p_0-p_1) \big)^{-2} \, \omega^{-\frac{3}{2}} \nonumber \\
												& \label{R221} = \frac{\sqrt{\pi}}{2^{\mu-2}} \, \| \tilde{u} \|_{L^{\infty}(p_1,p_2)} (p_0-p_1)^{\mu-3} \, \omega^{-\frac{3}{2}} \; . 
		\end{align}
	\end{itemize}
	To conclude, we give the expressions of $\tilde{K}_{\mu}^{(1)}(\omega,p_0,U), \tilde{H}^{(1)}(\omega,p_0,U)$ from \eqref{B11} and \eqref{B12}:
	\begin{align}
		& \label{H1mu} \bullet \quad \tilde{K}_{\mu}^{(1)}(\omega,p_0,U) := \frac{\Gamma(\mu)}{2^{\mu}} \, e^{i \frac{\pi \mu}{2}} \, e^{i \omega \psi(p_1)} \, \tilde{u}(p_1) \; , \\
		& \label{H112} \bullet \quad \tilde{H}^{(1)}(\omega,p_0,U) := \frac{\sqrt{\pi}}{2} \, e^{-i \frac{\pi}{4}} \, e^{i \omega c} \, \tilde{u}(p_0) \; .
	\end{align}
	And the constants $R_k^{(1)}(U)$ and the exponents $\alpha_k^1, \alpha_k^2$ are given by \eqref{R111}, \eqref{R211}, \eqref{R121} and \eqref{R221}.\\
	
	\noindent \textit{Study of $I^{(2)}(\omega,p_0)$}. Firstly we remark that $\psi'$ is negative for all $p \in [p_0,p_2]$. To apply Theorem \ref{THM1}, we make the phase increasing by using $p \mapsto -p$. We get
	\begin{equation*}
		I^{(2)}(\omega,p_0) = \int_{p_0}^{p_2} U(p) e^{i \omega \psi(p)} \, dp = \int_{\check{p}_2}^{\check{p}_0} \check{U}(p) e^{i \omega \breve{\psi}(p)} \, dp \; ,
	\end{equation*}
	where we put $\check{U}(p) := U(-p)$, $\check{\psi}(p) := \psi(-p)$, $\check{p}_0 := -p_0$ and $\check{p}_2 := -p_2$.\\
	Due to this substitution, $\check{\psi}$ is now an increasing function that satisfies Assumption \emph{(P$_{1,2,N}$)} on $[\check{p}_2, \check{p}_0]$, for all $N \geqslant 1$, and by hypothesis $\check{U}$ verifies \emph{(A$_{1,1,1}$)} on $[\check{p}_2, \check{p}_0]$. Furthermore we remark that $p_2$ is not a singularity of the amplitude and not a stationary point. From this observation, it follows that the use of a cutting-point is not necessary. Hence, in the notations of Theorem \ref{THM1}, we employ only the expansion of the integral $\tilde{I}^{(2)}(\omega,p_0)$, leading to more precise estimates. So we obtain from \eqref{I_2},
	\begin{equation*}
		\begin{aligned}
			I^{(2)}(\omega,p_0)	& = \phi^{(2)}(0,\omega,2,1) \, k_2(0) \, e^{i \omega \check{\psi}(\check{p}_0)} - \phi^{(2)}(s_2,\omega,2,1) \, k_2(s_2) \, e^{i \omega \check{\psi}(\check{p}_0)} \\
								& \qquad \qquad + e^{i \omega \check{\psi}(\check{p}_0)} \int_0^{s_2} \phi^{(2)}(s,\omega,2,1) (k_2)'(s) \, ds \; .
		\end{aligned}
	\end{equation*}
	Let us clarify the first terms by studying the function $\varphi_2$. The definition of $\varphi_2$ and the expression of $\psi$ yield
	\begin{equation*}
		\varphi_2(p) = \big( \check{\psi}(\check{p}_0) - \check{\psi}(p) \big)^{\frac{1}{2}} = (\check{p}_0 - p) \; ,
	\end{equation*}
	for all $p \in [\check{p}_2, \check{p}_0]$. It follows that $(\varphi_2)'(p) = (\varphi_2^{-1})'(s) = - 1$, and by the definition of $k_2$, we obtain
	\begin{equation*}
		k_2(s) = \check{U}\big(\varphi_2^{-1}(s) \big) (\varphi_2^{-1})'(s) = -\check{U}\big(\varphi_2^{-1}(s) \big) \; .
	\end{equation*}
	Since $\varphi_2(\check{p_0}) = 0$ and $\varphi_2(\check{p_2}) = s_2$, we have $k_2(0) = -\check{U}(\check{p_0}) = - U(p_0)$ and $k_2(s_2) = -U(p_2) = 0$, by the hypothesis on $U$.\\ Combining these two results with the expression of $\Theta^{(2)}(2,1)$ coming from Theorem \ref{THM1}, we obtain
	\begin{equation} \label{B22}
		\begin{aligned}
			I^{(2)}(\omega,p_0)	& = \frac{\sqrt{\pi}}{2} \, e^{-i \frac{\pi}{4}} \, e^{i \omega \psi(p_0)} \, \tilde{u}(p_0) \, (p_0 - p_1)^{\mu-1} \, \omega^{-\frac{1}{2}} \\
								& \qquad \qquad + e^{i \omega \psi(p_0)} \int_0^{s_2} \phi^{(2)}(s,\omega,2,1) (k_2)'(s) \, ds \; .
		\end{aligned}
	\end{equation}
	
	As in the preceding step, let us estimate the remainder term. Firstly, we bound the number $s_2$ as follows:
	\begin{equation} \label{s22}
		s_2 = \check{p}_0 - \check{p}_2 = p_2 - p_0 \leqslant p_2 - p_1 \; .
	\end{equation}
	Now we establish an estimate of the first derivative of $k_2$. By the definition of this function and by the fact that $(\varphi_2^{-1})'=-1$, we have
	\begin{equation*}
		(k_2)'(s) = (1-\mu) \left(\check{p}_1 - \varphi_2^{-1}(s) \right)^{\mu-2} \check{\tilde{u}}\left(\varphi_2^{-1}(s) \right) + \left(\check{p}_1 - \varphi_2^{-1}(s) \right)^{\mu-1} \Big(\check{\tilde{u}}\Big)' \left(\varphi_2^{-1}(s)\right) \; ,
	\end{equation*}
	where $\check{p}_1 := -p_1$. Since $\varphi_2^{-1}(s) \in [\check{p}_2,\check{p}_0]$, it follows
	\begin{equation*}
		\left\| (k_2)' \right\|_{L^{\infty}(0,s_2)} \leqslant \left( (1-\mu) (p_0 - p_1)^{\mu-2} + (p_0 - p_1)^{\mu-1} \right) \left\| \tilde{u} \right\|_{W^{1,\infty}(p_1,p_2)} \; .
	\end{equation*}
	Combine these considerations with the estimate of the remainder given in Theorem \ref{THM2} to obtain
	\begin{align}
		\left| 	R_1^{(2)}(\omega,p_0) \right|	& \leqslant \frac{L_{\gamma,2}}{1-\gamma} \, (p_2-p_1)^{1-\gamma} \left\| \tilde{u} \right\|_{W^{1,\infty}(p_1,p_2)} \nonumber \\
												& \label{Re12} \qquad \qquad \times \left( (1-\mu)(p_0-p_1)^{\mu-2} + (p_0-p_1)^{\mu-1} \right) \omega^{-\delta} \; ,
	\end{align}
	where $\gamma \in (0,1)$ and $\delta = \frac{1}{2} \, (\gamma + 1) \in \big(\frac{1}{2}, 1 \big)$.\\
	The inequality \eqref{Re12} furnishes the constants $R_k^{(2)}(U)$ and the exponents $\beta_k^1, \beta_k^2$. And the coefficient $\tilde{H}^{(2)}(\omega,p_0,U)$ is defined by
	\begin{equation} \label{H212}
		\bullet \quad \tilde{H}^{(2)}(\omega,p_0,U) := \frac{\sqrt{\pi}}{2} \, e^{-i \frac{\pi}{4}} \, e^{i \omega c} \, \tilde{u}(p_0) \; .
	\end{equation}
\end{proof}


\vspace{0.5cm}

In the next result, we shall substantially use the preceding asymptotic expansions. Since the blow-up comes from the term $(p_0-p_1)^{-\alpha}$, for a certain $\alpha > 0$, we suppose that $p_0 - p_1 = \omega^{-\varepsilon}$, where $\varepsilon > 0$. Under this hypothesis, it is possible to compare the blow-up and the decay. The aim of the second result is to find values of $\varepsilon$ so that the remainder estimates tend to $0$ with a better decay rate than the one of the first term.

\begin{thm4} \label{THM4}
	Under the assumptions of Theorem \ref{THM3}, let $\delta \in \big[ \frac{\mu+1}{2}, 1 \big)$, take $\varepsilon \in \big(0,\delta - \frac{1}{2}\big)$ and suppose that $p_0 := p_1 + \omega^{-\varepsilon}$ for $\omega > (p_2 - p_1)^{-\frac{1}{\varepsilon}}$. Then
	\begin{itemize}
		\item there exist $\tilde{K}_{\mu}^{(1)}(\omega,p_0,U), \tilde{H}^{(1)}(\omega,p_0,U) \in \C$ and a constant $\tilde{R}^{(1)}(U) \geqslant 0$ independent on $\omega$ and $p_0$, such that the following estimate holds:
	\begin{equation*}
			\Big| I^{(1)}(\omega,p_0) - \tilde{K}_{\mu}^{(1)}(\omega,p_0,U) \, \omega^{-\mu + \varepsilon \mu} -	\tilde{H}^{(1)}(\omega,p_0,U) \, \omega^{-\frac{1}{2} + \varepsilon(1-\mu)} \Big| \leqslant \tilde{R}^{(1)}(U) \, \omega^{-\alpha(\varepsilon,\delta)} \; ,
	\end{equation*}
	where $-\alpha(\varepsilon,\delta) < \min\big\{-\mu + \varepsilon \mu, -\frac{1}{2} + \varepsilon(1-\mu)\big\}$;
		\item there exist $\tilde{H}^{(2)}(\omega,p_0,U) \in \C$ and a constant $\tilde{R}^{(2)}(U) \geqslant 0$ independent on $\omega$ and $p_0$, such that the following estimate holds:
	\begin{equation*}
		\Big| I^{(2)}(\omega,p_0) -	 \tilde{H}^{(2)}(\omega,p_0,U) \, \omega^{-\frac{1}{2} + \varepsilon(1-\mu)} \Big| \leqslant \tilde{R}^{(2)}(U) \, \omega^{-\beta(\varepsilon,\delta)} \; ,
	\end{equation*}
	where $-\beta(\varepsilon,\delta) < \big(-\frac{1}{2} + \varepsilon(1-\mu)\big)$.
	\end{itemize}
\end{thm4}

\begin{proof}
	The proof is based on the result of Theorem \ref{THM3}. We replace $p_0 - p_1$ by $\omega^{-\varepsilon}$ in the previous estimates and we compare the decay rates of the first terms and the ones of the remainders. Further the hypothesis $\omega > (p_2 - p_1)^{-\frac{1}{\varepsilon}}$ implies $p_0 \in (p_1,p_2)$, and so $I^{(1)}(\omega,p_0)$ and $I^{(2)}(\omega,p_0)$ are well-defined.\\
	
	\noindent \textit{Case $j=1$}. The first terms become
	\begin{equation*}
		\begin{aligned}
			& \bullet \quad \tilde{K}_{\mu}^{(1)}(\omega,p_0,U) \, (p_0 - p_1)^{-\mu} \, \omega^{-\mu} = \tilde{K}_{\mu}^{(1)}(\omega,p_0,U) \, \omega^{-\mu + \varepsilon \mu} \; , \\
			& \bullet \quad \tilde{H}^{(1)}(\omega,p_0,U) \, (p_0 - p_1)^{\mu - 1} \, \omega^{-\frac{1}{2}} = \tilde{H}^{(1)}(\omega,p_0,U) \, \omega^{-\frac{1}{2} + \varepsilon(1-\mu)} \; ,
		\end{aligned}
	\end{equation*}
	where the two coefficients are defined in \eqref{H1mu} and \eqref{H112}. Let us make a comment about the decay rates of the first terms: by the definitions of $\tilde{K}_{\mu}^{(1)}(\omega,p_0,U)$ and $\tilde{H}^{(1)}(\omega,p_0,U)$ and by the hypothesis $\tilde{u}(p_1) \neq 0$, these two coefficients can be bounded from above and below by a positive constant when $\omega$ is sufficiently large. It follows that the coefficients do not influence the decay and so the two first terms behave like $\omega^{-\mu + \varepsilon \mu}$ and $\omega^{-\frac{1}{2} + \varepsilon(1-\mu)}$ when $\omega$ tends to infinity.\\
	Concerning the remainder terms, we get the new estimates by using \eqref{R111}, \eqref{R211}, \eqref{R121} and \eqref{R221}:
	\begin{align}
		& \label{R11} \bullet \; \left| R_1^{(1)}(\omega,p_0) \right| \leqslant \frac{2^{1-\mu}}{\mu} \, \| \tilde{u} \|_{W^{1,\infty}(p_1,p_2)} \Big( 2(2-\mu) \, \omega^{-1+\varepsilon(2-\mu)} + \omega^{-1+\varepsilon(1-\mu)} \Big) \; , \\
		& \label{R21} \bullet \; \left| R_2^{(1)}(\omega,p_0) \right| \leqslant \frac{1-\mu}{2^{\mu-2}} \, \| \tilde{u} \|_{L^{\infty}(p_1,p_2)} \, \omega^{-2 + \varepsilon(4-\mu)} \; , \\
		& \label{R12} \bullet \; \left| R_1^{(2)}(\omega,p_0) \right| \leqslant \frac{L_{\gamma, 2}}{1-\gamma} \, 2^{\gamma-\mu} \, \| \tilde{u} \|_{W^{1,\infty}(p_1,p_2)} \Big( 2(1-\mu) \, \omega^{-\delta + \varepsilon(1+\gamma-\mu)} + \omega^{-\delta + \varepsilon(\gamma-\mu)} \Big) \; , \\
		& \label{R22} \bullet \; \left| R_2^{(2)}(\omega,p_0) \right| \leqslant \frac{\sqrt{\pi}}{2^{\mu-2}} \, \| \tilde{u} \|_{L^{\infty}(p_1,p_2)} \, \omega^{-\frac{3}{2}+\varepsilon(3-\mu)} \; .
	\end{align}
	Note that if $\varepsilon \in \big(0, \delta - \frac{1}{2} \big)$ then $-\mu + \varepsilon \mu$ and $-\frac{1}{2} + \varepsilon (1-\mu)$ are both negative and so the first terms tend to $0$ like $\omega^{-\mu + \varepsilon \mu}$ and $\omega^{-\frac{1}{2} + \varepsilon (1-\mu)}$. Moreover
	\begin{equation*}
		-\frac{1}{2} + \varepsilon(1-\mu) \leqslant -\mu + \varepsilon \mu \qquad \Longleftrightarrow \qquad \mu \leqslant \frac{1}{2} \; .
	\end{equation*}
	Since the highest decay rate depends on the value of $\mu$, we shall distinguish the two cases.
	\begin{itemize}
	\item \textit{Case $\mu \leqslant \frac{1}{2}$.} In this situation, $\omega^{-\frac{1}{2} + \varepsilon(1-\mu)}$ decreases quicker than $\omega^{-\mu + \varepsilon \mu}$. We require the remainders to be asymptotically negligible as compared with $\omega^{-\frac{1}{2} + \varepsilon(1-\mu)}$ for $\omega$ tending to infinity. For this purpose, we show that the six exponents of $\omega$ that appear in the previous estimates are smaller than $-\frac{1}{2} + \varepsilon(1-\mu)$. \\ Estimates \eqref{R11}, \eqref{R21}, \eqref{R12} and \eqref{R22} lead to the following system of inequalities:
	\begin{subnumcases}{}
		\displaystyle -\frac{1}{2} + \varepsilon(1-\mu) > -1 + \varepsilon(2-\mu) \label{systineq2} \\
		\displaystyle -\frac{1}{2} + \varepsilon(1-\mu) > -1 + \varepsilon(1-\mu) \label{systineq3} \\
		\displaystyle -\frac{1}{2} + \varepsilon(1-\mu) > -2 + \varepsilon(4-\mu) \label{systineq4} \\
		\displaystyle -\frac{1}{2} + \varepsilon(1-\mu) > -\delta + \varepsilon(1+\gamma-\mu) \label{systineq6} \\
		\displaystyle -\frac{1}{2} + \varepsilon(1-\mu) > -\delta + \varepsilon(\gamma - \mu) \label{systineq5} \\
		\displaystyle -\frac{1}{2} + \varepsilon(1-\mu) > - \frac{3}{2} + \varepsilon(3-\mu) \label{systineq1}
	\end{subnumcases}
	Let us study each inequality.
	\begin{itemize}
		\item \emph{Inequalities \eqref{systineq2}, \eqref{systineq3}, \eqref{systineq4} and \eqref{systineq1}.} One can easily see that these inequalities are satisfied for $\varepsilon \in \big(0, \frac{1}{2} \big)$, which is true since $\varepsilon \in \big(0, \delta - \frac{1}{2} \big)$.
		\item \emph{Inequality \eqref{systineq6}.} This is equivalent to $\gamma^{-1} \big( \delta - \frac{1}{2} \big) > \varepsilon$. But $\gamma = 2 \delta - 1$ so \eqref{systineq6} is equivalent to $\varepsilon < \frac{1}{2}$, that is true by hypothesis.
		\item \emph{Inequality \eqref{systineq5}.} This inequality is equivalent to $ \delta - \frac{1}{2} + \varepsilon(1-\gamma) > 0$. The left-hand side of the last inequality is an increasing function with respect to $\varepsilon > 0$ since $\gamma < 1$. Furthermore, it is assumed that $\delta \in \big[ \frac{\mu+1}{2},1\big)$, hence $\delta > \frac{1}{2}$; so \eqref{systineq5} is verified for $\varepsilon > 0$.
	\end{itemize}
	\item \textit{Case $\mu \geqslant \frac{1}{2}$.} Here the quickest decay rate of the first terms is $\omega^{-\mu + \varepsilon \mu}$. We obtain a new system of inequalities:
	\begin{subnumcases}{}
		\displaystyle -\mu + \varepsilon \mu > -1 + \varepsilon(2-\mu) \label{systineq8} \\
		\displaystyle -\mu + \varepsilon \mu > -1 + \varepsilon(1-\mu) \label{systineq9} \\
		\displaystyle -\mu + \varepsilon \mu > -2 + \varepsilon(4-\mu) \label{systineq10} \\
		\displaystyle -\mu + \varepsilon \mu > -\delta + \varepsilon(1+\gamma-\mu) \label{systineq12} \\
		\displaystyle -\mu + \varepsilon \mu > -\delta + \varepsilon(\gamma - \mu) \label{systineq11} \\
		\displaystyle -\mu + \varepsilon \mu > - \frac{3}{2} + \varepsilon(3-\mu) \label{systineq7}
	\end{subnumcases}
	We study these inequalities:
	\begin{itemize}
		\item \emph{Inequalities \eqref{systineq8}, \eqref{systineq10} and \eqref{systineq7}.} By a quick calculation, one can show that all are satisfied for $0 < \varepsilon < \frac{1}{2}$.
		\item \emph{Inequality \eqref{systineq9}.} This inequality is equivalent to $1-\mu + \varepsilon(2\mu-1) > 0$. Since $\frac{1}{2} \leqslant \mu < 1$, the left-hand side is an increasing function with respect to $\varepsilon$ which is strictly positive for all $\varepsilon > 0$.
		\item \emph{Inequality \eqref{systineq12}.} This one is equivalent to $\delta - \mu + \varepsilon(2\mu -1 -\gamma) > 0$. Note that if $\delta > \mu$ then $2\mu - 1 - \gamma < 0$ by using the relation $\delta = \frac{1}{2}(1+\gamma)$. So $\varepsilon \longmapsto \delta - \mu + \varepsilon(2\mu -1 -\gamma)$ is a decreasing function and it is positive if and only if $\varepsilon < \frac{\delta-\mu}{\gamma+1-2\mu} = \frac{1}{2}$. But the assumption on $\delta$ implies $\delta > \mu$, so \eqref{systineq12} is true for $\varepsilon < \frac{1}{2}$.
		\item \emph{Inequality \eqref{systineq11}.} This is equivalent to $\delta - \mu + \varepsilon(2 \mu - \gamma) > 0$. Since $\gamma \in (0,1)$ and $\frac{1}{2} \leqslant \mu$, the left-hand side is an increasing function with respect to $\varepsilon$. Moreover we have $\delta > \mu$, that assures the positivity of the left-hand side for all $\varepsilon > 0$.
	\end{itemize}
	\end{itemize}
	To finish the case $j=1$, we choose $\alpha(\varepsilon,\delta)$ such that $\omega^{-\alpha(\varepsilon,\delta)}$ is the lowest decay rate of the remainder. This implies
	\begin{equation*}
	 	\sum_{k = 1}^6 \, R_k^{(1)}(U) \, \omega^{-\alpha_k^2 + \varepsilon \alpha_k^1} = \omega^{-\alpha(\varepsilon,\delta)} \sum_{k = 1}^6 \, R_k^{(1)}(U) \, \omega^{-\alpha_k^2 + \varepsilon \alpha_k^1 + \alpha(\varepsilon,\delta)} \leqslant \tilde{R}^{(1)}(U) \, \omega^{-\alpha(\varepsilon,\delta)} \; .
	\end{equation*}	
	since each $\omega^{-\alpha_k^2 + \varepsilon \alpha_k^1 + \alpha(\varepsilon,\delta)}$ is bounded when $\omega > (p_2 - p_1)^{-\frac{1}{\varepsilon}}$.\\
	
	\noindent \textit{Case $j=2$}. In the same manner, we obtain
	\begin{equation*}
		\bullet \quad \tilde{H}^{(2)}(\omega,p_0,U) \, (p_0-p_1)^{\mu-1} \, \omega^{-\frac{1}{2}} = \tilde{H}^{(2)}(\omega,p_0,U) \, \omega^{-\frac{1}{2} + \varepsilon(1-\mu)} \; ,
	\end{equation*}
	by putting $p_0 - p_1 = \omega^{-\varepsilon}$ where the coefficient is defined in \eqref{H212}. As above, we observe that this term behaves like $\omega^{-\frac{1}{2} + \varepsilon(1-\mu)}$ when $\omega$ tends to infinity. From \eqref{Re12}, we obtain the new estimate:
	\begin{equation} \label{R122}
		\bullet \; \left| R_1^{(2)}(\omega,p_0) \right| \leqslant \frac{L_{\gamma,2}}{1-\gamma} \, (p_2-p_1)^{1-\gamma} \left\| \tilde{u} \right\|_{W^{1,\infty}(p_1,p_2)} \Big( (1-\mu) \, \omega^{-\delta + \varepsilon(2-\mu)} + \omega^{-\delta + \varepsilon (1-\mu)} \Big) \; .
	\end{equation}
	Let us compare the decay rates and remark that we do not have to distinguish different cases because there is only one first term here. Estimate \eqref{R122} leads to the following system of inequalities:
	\begin{subnumcases}{}
		\displaystyle -\frac{1}{2} + \varepsilon(1-\mu) > -\delta + \varepsilon(1-\mu) \label{systineq14} \\
		\displaystyle -\frac{1}{2} + \varepsilon(1-\mu) > -\delta + \varepsilon(2-\mu) \label{systineq15}
	\end{subnumcases}
	\begin{itemize}
		\item \textit{Inequality \eqref{systineq14}.} Clearly this is equivalent to $\delta > \frac{1}{2}$, which is satisfied.
		\item \textit{Inequality \eqref{systineq15}.} This inequality is equivalent to $\delta - \frac{1}{2} - \varepsilon > 0$, which is true by the hypothesis on $\varepsilon$.
	\end{itemize}
	To complete the proof, $\beta(\varepsilon,\delta)$ is chosen in such a way that $\omega^{-\beta(\varepsilon,\delta)}$ is the lowest decay rate of the remainder. It follows
	\begin{equation*}
	 	\sum_{k = 1}^2 \, R_k^{(2)}(U) \, \omega^{-\beta_k^2 + \varepsilon \beta_k^1} = \omega^{-\beta(\varepsilon,\delta)} \sum_{k = 1}^2 \, R_k^{(2)}(U) \, \omega^{-\beta_k^2 + \varepsilon \beta_k^1 + \beta(\varepsilon,\delta)} \leqslant \tilde{R}^{(2)}(U) \, \omega^{-\beta(\varepsilon,\delta)} \; ,
	\end{equation*}	
	by the boundedness of $\omega^{-\beta_k^2 + \varepsilon \beta_k^1 + \beta(\varepsilon,\delta)}$ for $\omega > (p_2 - p_1)^{-\frac{1}{\varepsilon}}$.\\
	The result is proved.
\end{proof}

\begin{rem3}
	\em Note that if $\varepsilon = \frac{1}{2}$ then the decay rates of the expansion and the remainder are equal according to the proof. Moreover if $\varepsilon$ tends to $\frac{1}{2}$ then $\delta$ have to tend to $1$ and so $\gamma$ tends to $1$ as well. In this situation, the estimates \eqref{R12} and \eqref{R122} blow up, implying that the constants $\tilde{R}^{(1)}(U), \tilde{R}^{(2)}(U)$ tend to infinity when $\varepsilon$ tends to $\frac{1}{2}$.
\end{rem3}

\section{The free Schrödinger equation on the line: low decay in space-time regions approaching the critical speed associated with a frequency domain singularity of the initial condition}

\hspace{2.5ex} In this section, we are interested in the time-asymptotic behaviour of the solution of the free Schrödinger equation in one dimension, where the Fourier transform of the initial data is supported on $[p_1,p_2]$ and has a singularity at $p_1$. In \cite{article1}, we furnished asymptotic expansions outside narrow cones centred in the space-time directions given by the endpoints of the frequency band, namely $p_1$ and $p_2$. A blow-up occurs when the boundaries of the cones tend to the critical directions and the employed method did not allow us to derive optimal decay rates in larger regions.

The aim of this section is to replace the uncontrollable cone generated by the singular frequency by a smaller region, using the results of section 2. To do so, we write the solution of the free Schrödinger equation as an oscillatory integral with respect to time. We observe that the phase function has a unique stationary point $p_0 := \frac{x}{2t}$, which furnishes a family of oscillatory integrals indexed by the parameter $p_0$. Consequently, the results of the preceding section are substantially used.

In the first result, Theorem \ref{THM3} is employed to obtain explicit estimates of the integral given by the solution. Since the blow-up depends only on the distance between the singularity $p_1$ and the stationary point $p_0$, we suppose that this distance is bounded from below by $t^{-\varepsilon}$, with $\varepsilon > 0$, leading to a space-time region $\mathfrak{R}_{\varepsilon}$. According to the previous section, if the parameter $\varepsilon$ is sufficiently small, then there is no blow-up and we can uniformly estimate the solution in the above mentioned space-time region, which  is asymptotically larger than any space-time cone included between the directions $p_1$ and $p_2$.

The second result is devoted to the optimality of the previous estimates. In the region $\mathfrak{R}_{\varepsilon}$, we expect that the decay rate will be slow in parts which are close to the critical direction given by $p_1$, where the influence of the singularity is the strongest. So we use Theorem \ref{THM4} to provide asymptotic expansions of the solution on the space-time curve $\mathfrak{S}_{\varepsilon}$, the left boundary of the region $\mathfrak{R}_{\varepsilon}$, and we show that the decay rates obtained in the preceding result are attained on $\mathfrak{S}_{\varepsilon}$, proving the optimality.\\

First of all we describe the curve $\mathfrak{G}_{\varepsilon}$ in the following proposition. The point $iii)$ affirms in particular that the region $\mathfrak{R}_{\varepsilon}$ is included between the directions $p_1$ and $p_2$ and is asymptotically larger than any space-time cone.

\begin{prop1} \label{prop1}
	Let $\varepsilon > 0$ and $p_1 \in \R$. Define the curve $\mathfrak{G}_{\varepsilon} \subset \R_+^* \times \R$ as follows:
	\begin{equation*}
		\mathfrak{G}_{\varepsilon} := \left\{ (t,x) \in \R_+^* \times \R \: \Big| \: \frac{x}{2t} - p_1 = t^{-\varepsilon} \right\} \; .
	\end{equation*}
	Then the following assertions are true:
	\begin{enumerate}
		\item $\forall \, (t,x) \in \mathfrak{G}_{\varepsilon} \qquad	x = 2 p_1 \, t + t^{1-\varepsilon} \; .$
		\item $\forall \, (t,x) \in \mathfrak{G}_{\varepsilon} \qquad x > 2 p_1 \, t \; ;$
		\item for all $p > p_1$ and for all $(t,x) \in \mathfrak{G}_{\varepsilon}$ where $t \geqslant T_{p} := \big(2(p-p_1)\big)^{-\frac{1}{\varepsilon}} > 0$, we have
		\begin{equation*}
			x \leqslant 2 p \, t \; ;
		\end{equation*}
	\end{enumerate}
\end{prop1}

\begin{proof}
	Simple consequences of the definition of $\mathfrak{G}_{\varepsilon}$.
\end{proof}

\vspace{0.5cm}

Now we introduce the assumption concerning the initial data.\\

\noindent \textbf{Condition} ($C_{p_1,p_2,\mu}$): Fix $\mu \in (0,1)$ and let $p_1, p_2$ be two real numbers such that
	\begin{equation*}
		- \infty < p_1 < p_2 < +\infty \; .
	\end{equation*}
	$u_0$ satisfies condition ($C_{p_1,p_2,\mu}$) if and only if $\tf u_0 \equiv 0$ on $\R \setminus [p_1,p_2]$ and $U := \tf u_0$ verifies Assumption (A$_{\mu,1,1}$) on $[p_1,p_2]$, with $\tf u_0(p_2) = 0$.\\


\begin{cor2}
	Suppose that $u_0$ satisfies Condition ($C_{p_1,p_2,\mu}$). Fix $\delta \in \big[ \frac{\mu+1}{2},1\big)$ and $\varepsilon \in \big(0, \delta - \frac{1}{2} \big)$. Then there exist three constants $C_{\mu > 1/2}(u_0), C_{\mu = 1/2}(u_0), C_{\mu < 1/2}(u_0) > 0$ such that for all $(t,x) \in \mathfrak{R}_{\varepsilon}$, the following estimates hold:
	\begin{itemize}
	\item Case $\mu > \frac{1}{2}$:
	\begin{equation*}
		\big| u(t,x) \big| \leqslant C_{\mu > 1/2}(u_0) \, t^{-\frac{1}{2} + \varepsilon(1-\mu)} \; ;
	\end{equation*}
	\item Case $\mu = \frac{1}{2}$:
	\begin{equation*}
		\big| u(t,x) \big| \leqslant C_{\mu = 1/2}(u_0) \,  t^{-\frac{1}{2} + \frac{\varepsilon}{2}} \; ;
	\end{equation*}
	\item Case $\mu < \frac{1}{2}$:	
	\begin{equation*}
		\big| u(t,x) \big| \leqslant C_{\mu < 1/2}(u_0) \, t^{-\mu + \varepsilon \mu} \; .
	\end{equation*}
	\end{itemize}
	The region $\mathfrak{R}_{\varepsilon}$ is described as follows:
	\begin{equation*}
		(t,x) \in \mathfrak{R}_{\varepsilon} \qquad \Longleftrightarrow \qquad \left\{ \begin{array}{rl}
			& \displaystyle \frac{x}{2t} - p_1 \geqslant t^{-\varepsilon} \; ,\\
			& \displaystyle x < 2 p_2 \, t \; , \\
			& \displaystyle t > T_{p_2} \; .
		\end{array} \right.
	\end{equation*}
\end{cor2}

\begin{rem4}
	\em The definition of $\mathfrak{R}_{\varepsilon}$ shows that $\mathfrak{G}_{\varepsilon}$ is the left boundary while the direction $p_2$ is the right boundary. Further the hypothesis $t > T_{p_2}$ assures that $\mathfrak{R}_{\varepsilon}$ is non-empty.
\end{rem4}

\begin{proof}
	We divide the proof in three steps, inspired by \cite[section 2]{article1}.\\
	
	\noindent \textit{First step: Rewriting.} We recall the expression of the solution of the free Schrödinger equation:
	\begin{equation*}
		\forall \, (t,x) \in [0,+\infty) \times \R \qquad u(t,x) = \frac{1}{2\pi} \int_{p_1}^{p_2} \tf u_0(p) e^{-i t p^2 + i x p} dp \; .
	\end{equation*}
	We factorize the phase function by $t$ and we obtain
	\begin{equation*}
		\forall \, (t,x) \in [0,+\infty) \times \R \qquad u(t,x) = \frac{1}{2\pi} \int_{p_1}^{p_2} \tf u_0(p) e^{i t \Psi(p,t,x)} dp \; ,
	\end{equation*}
	where $\Psi(p,t,x) := - p^2 + \frac{x}{t} p = -\big( p - \frac{x}{2t} \big)^2 + \frac{x^2}{4 t^2}$.\\
	
	\noindent \textit{Second step: Splitting of the integral.} Let us define the amplitude $U$ and the phase $\psi$ by
	\begin{equation*}
	\left\{ \begin{array}{rl}
			& \displaystyle \forall \, p \in (p_1,p_2] \qquad U(p) := \tf u_0(p) = (p-p_1)^{\mu-1} \tilde{u}(p) \; , \\
			& \displaystyle \forall \, p \in [p_1,p_2] \qquad \psi(p) := \Psi(p,t,x) \; .
	\end{array} \right.
	\end{equation*}
	Then the first derivative $\psi'$ of the phase is given by
	\begin{equation*} \label{psi}
		\forall \, p \in [p_1,p_2] \qquad \psi'(p) = - 2 \Big(p - \frac{x}{2 t}\Big) =: 2 \big(p_0(t,x) -p\big) \; ,
	\end{equation*}
	where we put $p_0(t,x) := \frac{x}{2t}$, which is the only stationary point of the phase.\\
	The hypothesis $(t,x) \in \mathfrak{R}_{\varepsilon}$ implies that $p_0 := p_0(t,x) < p_2$ and $p_0 - p_1 \geqslant t^{-\varepsilon}$. In such a situation, $p_0 \in (p_1,p_2)$ and we can write
	\begin{equation*}
		u(t,x) = \frac{1}{2\pi} \bigg( \int_{p_1}^{p_0} U(p) e^{i t \psi(p)} dp \; + \; \int_{p_0}^{p_2} U(p) e^{i t \psi(p)} dp \bigg) =: \frac{1}{2\pi} \big( I^{(1)}(t,p_0) + I^{(2)}(t,p_0) \big) \; .
	\end{equation*}
	Before applying Theorem \ref{THM3}, we assume that $\mu > \frac{1}{2}$; the cases $\mu = \frac{1}{2}$ and $\mu < \frac{1}{2}$ can be treated in a similar way.\\
	
	\noindent \textit{Third step: Application of Theorem \ref{THM3}}. By putting $p_0 := \frac{x}{2t}$, $c := \frac{x^2}{4 t^2}$ and by choosing $t$ as the large parameter, the hypothesis of Theorem \ref{THM3} are satisfied and the following estimates hold:
	\begin{equation*}
		\begin{aligned}
			& \bullet \quad \left| I^{(1)}(t,p_0) \right| \leqslant \left| \tilde{H}^{(1)}(t, p_0, U) \right| (p_0-p_1)^{\mu-1} \, t^{-\frac{1}{2}} \\
			& \qquad + \left| \tilde{K}_{\mu}^{(1)}(t, p_0, U) \right| (p_0-p_1)^{-\mu} \, t^{-\mu} + \sum_{k=1}^6 R_k^{(1)}(U) (p_0-p_1)^{-\alpha_k^1} \, t^{-\alpha_k^2} \; , \\
			& \bullet \quad \left| I^{(2)}(t,p_0) \right| \leqslant \left| \tilde{H}^{(2)}(t, p_0, U) \right| (p_0-p_1)^{\mu-1} \, t^{-\frac{1}{2}} + \sum_{k=1}^2 R_k^{(2)}(U) (p_0-p_1)^{-\beta_k^1} \, t^{-\beta_k^2} \; .
		\end{aligned}
	\end{equation*}
	Combining the hypothesis $p_0 - p_1 \geqslant t^{-\varepsilon}$ and the fact that $\alpha_k^1, \beta_k^1 \geqslant 0$ leads to
	\begin{align}
		& \bullet \quad \left| I^{(1)}(t,p_0) \right| \leqslant \left| \tilde{H}^{(1)}(t, p_0, U) \right| \, t^{-\frac{1}{2} + \varepsilon(1-\mu)} \nonumber \\
		& \label{EstI1} \qquad + \left| \tilde{K}_{\mu}^{(1)}(t, p_0, U) \right| \, t^{-\mu + \varepsilon \mu} + \sum_{k=1}^6 R_k^{(1)}(U) \, t^{-\alpha_k^2 + \varepsilon \alpha_k^1} \; , \\
		& \label{EstI2} \bullet \quad \left| I^{(2)}(t,p_0) \right| \leqslant \left| \tilde{H}^{(2)}(t, p_0, U) \right| \, t^{-\frac{1}{2} + \varepsilon(1-\mu)} + \sum_{k=1}^2 R_k^{(2)}(U) \, t^{-\beta_k^2 + \varepsilon \beta_k^1} \; .
	\end{align}
	Moreover the definitions of $\tilde{H}^{(1)}(t,p_0,U) = \tilde{H}^{(2)}(t,p_0,U)$ and $\tilde{K}_{\mu}^{(1)}(t,p_0,U)$ give
	\begin{equation*}
		\begin{aligned}
			& \bullet \quad \left| \tilde{H}^{(j)}(t, p_0, U) \right| \leqslant \frac{\sqrt{\pi}}{2} \, \| \tilde{u} \|_{L^{\infty}(p_1,p_2)} \; , \quad j = 1,2 \; , \\
			& \bullet \quad \left| \tilde{K}_{\mu}^{(1)}(t, p_0, U) \right| \leqslant \frac{\Gamma(\mu)}{2^{\mu}} \, \| \tilde{u} \|_{L^{\infty}(p_1,p_2)} \; .
		\end{aligned}
	\end{equation*}
	Theorem \ref{THM4} affirms that $t^{-\frac{1}{2} + \varepsilon(1-\mu)}$ is the lowest decay rate if $\delta \in \big[ \frac{\mu + 1}{2}, 1 \big)$ and $\varepsilon \in \big(0, \delta - \frac{1}{2}\big)$. Hence we can factorize \eqref{EstI1} and \eqref{EstI2} by $t^{-\frac{1}{2} + \varepsilon(1-\mu)}$ and we use the fact that $t > T_{p_2}$ to obtain
	\begin{align}
		& \label{I1} \bullet \quad \left| I^{(1)}(t,p_0) \right| \leqslant C_{\mu > 1/2}^{(1)}(u_0) \, t^{-\frac{1}{2} + \varepsilon(1-\mu)} \; , \\
		& \label{I2} \bullet \quad \left| I^{(2)}(t,p_0) \right| \leqslant C_{\mu > 1/2}^{(2)}(u_0) \, t^{-\frac{1}{2} + \varepsilon(1-\mu)} \; ,
	\end{align}
	for certain constants $C_{\mu > 1/2}^{(j)}(u_0) \geqslant 0$ coming from \eqref{EstI1} and \eqref{EstI2}. Finally we obtain the desired estimate for $u(t,x)$ for all $(t,x) \in \mathfrak{R}_{\varepsilon}$ by adding \eqref{I1} and \eqref{I2}:
	\begin{equation*}
		\begin{aligned}
			\big| u(t,x) \big|	& \leqslant \frac{1}{2\pi} \left| I^{(1)}(t,p_0) \right| + \frac{1}{2 \pi} \left| I^{(2)}(t,p_0) \right| \\
								& \leqslant \frac{C_{\mu > 1/2}^{(1)}(u_0)}{2\pi} \, t^{-\frac{1}{2} + \varepsilon(1-\mu)} + \frac{C_{\mu > 1/2}^{(2)}(u_0)}{2 \pi} \, t^{-\frac{1}{2} + \varepsilon(1-\mu)} \\
								& =: C_{\mu > 1/2}(u_0) \, t^{-\frac{1}{2} + \varepsilon(1-\mu)} \; .
		\end{aligned}
	\end{equation*}
\end{proof}

\vspace{0.1cm}

\begin{thm5} \label{THM5}
	Suppose that $u_0$ satisfies Condition ($C_{p_1,p_2,\mu}$). Fix $\delta \in \big[ \frac{\mu+1}{2},1\big)$ and $\varepsilon \in \big(0, \delta - \frac{1}{2} \big)$. Then for all $(t,x) \in \mathfrak{G}_{\varepsilon}$ where it is assumed $t > T_{p_2}$, there exist $K_{\mu}(t,u_0), H(t,u_0) \in \C$ uniformly bounded by $R_K(u_0), R_H(u_0) \geqslant 0$ respectively, satisfying
	\begin{itemize}
	\item Case $\mu > \frac{1}{2}$:	
	\begin{equation*}
		\Big| u(t,x) - H(t,u_0) \, t^{-\frac{1}{2} + \varepsilon(1-\mu)} \, \Big| \leqslant R_K(u_0) \, t^{-\mu + \varepsilon \mu} + R^{(1)}(u_0) \, t^{-\alpha(\varepsilon,\delta)} + R^{(2)}(u_0) \, t^{-\beta(\varepsilon,\delta)} \; ,
	\end{equation*}
	where $-\alpha(\varepsilon,\delta), -\beta(\varepsilon,\delta) < -\frac{1}{2} + \varepsilon (1-\mu)$.
	\item Case $\mu = \frac{1}{2}$:
	\begin{equation*}
		\Big| u(t,x) - \big(K_{\mu}(t,u_0) + H(t,u_0) \big) \, t^{-\frac{1}{2} + \frac{\varepsilon}{2}} \, \Big| \leqslant R^{(1)}(u_0) \, t^{-\alpha(\varepsilon,\delta)} + R^{(2)}(u_0) \, t^{-\beta(\varepsilon,\delta)} \; ,
	\end{equation*}
	where $-\alpha(\varepsilon,\delta), -\beta(\varepsilon,\delta) < -\frac{1}{2}+ \frac{\varepsilon}{2}$.
	\item Case $\mu < \frac{1}{2}$:	
	\begin{equation*}
		\Big| u(t,x) - K_{\mu}(t,u_0) \, t^{-\mu + \varepsilon \mu} \, \Big| \leqslant R_H(u_0) \, t^{-\frac{1}{2} + \varepsilon(1-\mu)} + R^{(1)}(u_0) \, t^{-\alpha(\varepsilon,\delta)} + R^{(2)}(u_0) \, t^{-\beta(\varepsilon,\delta)} \; ,
	\end{equation*}
	where $-\alpha(\varepsilon,\delta), -\beta(\varepsilon,\delta) < -\mu + \varepsilon \mu$.
	\end{itemize}
	$R^{(1)}(u_0)$, $R^{(2)}(u_0) \geqslant 0$ are constants independent on $t,x$.
\end{thm5}

\begin{proof}
	As above, we study only the case $\mu > \frac{1}{2}$; the others cases can be treated in a similar way.\\
	We shall use the rewriting of the solution coming from the two first steps of the preceding proof:
	\begin{equation*}
		u(t,x) = \frac{1}{2\pi} \int_{p_1}^{p_2} U(p) e^{i t \psi(p)} dp \; .
	\end{equation*}
	We recall the definition of $p_0 := \frac{x}{2t}$. $(t,x) \in \mathfrak{G}_{\varepsilon}$ is equivalent to $p_0 - p_1 = t^{-\varepsilon}$, and by Proposition \ref{prop1}, if $t > T_{p_2}$ then $2 p_1 \, t < x < 2 p_2 \, t$, that is to say $p_0 \in (p_1,p_2)$. Hence one can split the integral at the point $p_0$ and one obtains the integrals $I^{(1)}(t, p_0)$ and $I^{(2)}(t,p_0)$. $U$ and $\psi$ satisfy the hypothesis of Theorem \ref{THM4}, which is applicable with $t$ the large parameter, and we obtain the following estimates
	\begin{align*}
		& \bullet \quad \Big| I^{(1)}(t,p_0) - \tilde{H}^{(1)}(t,p_0,U) \, t^{-\frac{1}{2} + \varepsilon(1-\mu)} \Big| \leqslant \big| \tilde{K}_{\mu}^{(1)}(t,p_0,U) \big| \, t^{-\mu + \varepsilon \mu} + \tilde{R}^{(1)}(U) \, t^{-\alpha(\varepsilon,\delta)} \; , \\
		& \bullet \quad \Big| I^{(2)}(t,p_0) -	 \tilde{H}^{(2)}(t,p_0,U) \, t^{-\frac{1}{2} + \varepsilon(1-\mu)} \Big| \leqslant \tilde{R}^{(2)}(U) \, t^{-\beta(\varepsilon,\delta)} \; .
	\end{align*}
	Let us compute the coefficients. We have
	\begin{equation*}
		\begin{aligned}
			\tilde{H}^{(1)}(t,p_0,U) & =  \frac{\sqrt{\pi}}{2} \, e^{-i \frac{\pi}{4}} \, e^{i t (p_1 + t^{-\varepsilon})} \, \tilde{u}(p_1+ t^{-\varepsilon}) \\
									& =: \pi \, H(t,u_0) \; ,
		\end{aligned}
	\end{equation*}
	where $H(t,u_0)$ is uniformly bounded by $R_H(u_0) := \frac{1}{2\sqrt{\pi}} \| \tilde{u} \|_{L^{\infty}(p_1,p_2)}$. We carry out the same computation for $\tilde{H}^{(2)}(t, p_0, U)$ and we obtain the same result, namely
	\begin{equation*}
		\tilde{H}^{(2)}(t, p_0, U) := \pi \, H(t,u_0) \; .
	\end{equation*}
	Then the expression of $\tilde{K}_{\mu}^{(1)}(t, p_0, U)$ from Theorem \ref{THM3} (see \eqref{H1mu}) yields
	\begin{equation*}
		\left| \tilde{K}_{\mu}^{(1)}(t, p_0, U) \right|\leqslant \frac{\Gamma(\mu)}{2^{\mu}} \, \| \tilde{u} \|_{L^{\infty}(p_1,p_2)}  =: \tilde{R}_{K}(u_0) \; .
	\end{equation*}
	 To conclude, we set $R_K(u_0) := \frac{1}{2 \pi} \tilde{R}_{K}(u_0)$, $R^{(1)}(u_0) := \frac{1}{2 \pi} \tilde{R}^{(1)}(U)$, $R^{(2)}(u_0) := \frac{1}{2 \pi} \tilde{R}^{(2)}(U)$ and we obtain the final estimate
	\begin{equation*}
		\begin{aligned}
			\Big| u(t,x)	& - H(t,u_0) \, t^{-\frac{1}{2} + \varepsilon(1-\mu)} \Big|	\\
							& \leqslant \frac{1}{2\pi} \, \Big| I^{(1)}(t,p_0) - \tilde{H}^{(1)}(t, p_0, U) \, t^{-\frac{1}{2} + \varepsilon(1-\mu)} \Big| \\
							& \qquad \qquad + \frac{1}{2 \pi} \, \Big| I^{(2)}(t,p_0) - \tilde{H}^{(2)}(t, p_0, U) \, t^{-\frac{1}{2} + \varepsilon(1-\mu)} \Big| \\
							& \leqslant \frac{1}{2 \pi} \, \left| \tilde{K}_{\mu}^{(1)}(t, p_0, U) \right| \, t^{-\mu + \varepsilon \mu} + \frac{1}{2 \pi} \, \tilde{R}^{(1)}(U) \, t^{-\alpha(\varepsilon,\delta)} + \frac{1}{2 \pi} \, \tilde{R}^{(2)}(U) \, t^{-\beta(\varepsilon,\delta)} \\
							& \leqslant R_K(u_0) \, t^{-\mu + \varepsilon \mu} + R^{(1)}(u_0) \, t^{-\alpha(\varepsilon,\delta)} + R^{(2)}(u_0) \, t^{-\beta(\varepsilon,\delta)} \; .
		\end{aligned}
	\end{equation*}	
\end{proof}

\section{Appendix: technical results from part I}

\hspace{2.5ex} In this section, we formulate some results which are useful in section 1. Their statements are very similar to the ones from \cite{article1}, the only difference comes from the definition of the point $s_j$. But this does not affect the results and the proofs.

\begin{prop2} \label{PROP2}
	Fix $s > 0$ and $\rho_j \geqslant 1$. Let $\Lambda^{(j)}(s)$ be the curve of the complex plane given in Definition \ref{DEF1}. Then we have
	\begin{equation*}
		\forall \, z = s + t e^{(-1)^{j+1} i \frac{\pi}{2 \rho_j}} \in \Lambda^{(j)}(s) \qquad \left| \, e^{(-1)^{j+1} i \omega z^{\rho_j}} \right| \leqslant e^{- \omega t^{\rho_j}} \; .
	\end{equation*}
\end{prop2}

\begin{prop3} \label{PROP3}
	Let $\psi : [p_1,p_2] \longrightarrow \R$ be a function which satisfies Assumption \emph{(P$_{\rho_1,\rho_2,N}$)}. Consider the function $\varphi_j : I_j \longrightarrow \R$ introduced in Definition \ref{DEF1}. Then $\varphi_j$ is a $\mathcal{C}^{N+1}$-diffeomorphism between $I_j$ and $[0,s_j]$.
\end{prop3}

\begin{prop4} \label{PROP4}
	Let $U : (p_1,p_2) \longrightarrow \C$ be a function which satisfies Assumption \emph{(A$_{\mu_1,\mu_2,N}$)}. Consider the function $k_j : (0,s_j] \longrightarrow \C$ introduced in Definition \ref{DEF1}. Then $k_j$ can be extended to the interval $[0,s_j]$ and $k_j \in \mathcal{C}^{N}\big([0,s_j] \big)$.
\end{prop4}

\begin{cor1} \label{COR1}
	Fix $s_j > 0$, $\rho_j \geqslant 1$ and $0< \mu_j \leqslant 1$. For any $\omega >0$, the function $\phi^{(j)}(.,\omega,\rho_j,\mu_j) : \, (0,s_j] \longrightarrow \C$ defined by 
	\begin{equation*}
		\phi^{(j)}(s,\omega,\rho_j,\mu_j) = -\int_{\Lambda^{(j)}(s)} z^{\mu_j-1} e^{(-1)^{j+1} i \omega z^{\rho_j}} \, dz
	\end{equation*}
	is a primitive of the function $s \in (0,s_j] \longmapsto s^{\mu_j-1} e^{(-1)^{j+1} i \omega s^{\rho_j}} \in \C$.
\end{cor1}

\end{document}